\documentclass[11pt,a4paper]{article}
\usepackage[T1]{fontenc}
\usepackage{lmodern,amsmath,amssymb,amsthm,mathtools}
\usepackage[margin=24mm]{geometry}
\usepackage{microtype,xcolor}
\usepackage[colorlinks=true,allcolors=blue!55!black]{hyperref}
\hypersetup{pdftitle={Abelian covers and the index of the centre in a group},pdfauthor={Guillaume Lecomte}}

\theoremstyle{plain}
\newtheorem{theorem}{Theorem}[section]
\newtheorem{lemma}[theorem]{Lemma}
\newtheorem{proposition}[theorem]{Proposition}
\newtheorem{corollary}[theorem]{Corollary}
\newcommand{\casecell}[1]{{\raggedright #1\par}}

\newcommand{\F}{\mathbf F}
\newcommand{\Z}{\mathrm Z}
\DeclareMathOperator{\Tr}{Tr}
\newcommand{\doi}[1]{\href{https://doi.org/#1}{\nolinkurl{doi:#1}}}
\title{Abelian covers and the index of the centre in a group}
\author{Guillaume Lecomte\thanks{Independent researcher.
Correspondence: \href{mailto:guillaume.lecomteexed@edu.executive.em-lyon.com}{\nolinkurl{guillaume.lecomteexed@edu.executive.em-lyon.com}}.
ORCID: \href{https://orcid.org/0009-0003-3532-3379}{0009-0003-3532-3379}.}}
\date{}

\begin{document}
\maketitle

\begin{abstract}
Let \(h(n)\) be the largest number of abelian subgroups needed
to cover a group in which every set of more than \(n\) elements
contains two distinct elements that commute, and let \(\zeta(n)\)
be the largest index of the centre under the same constraint.
Erd\H{o}s asked for the growth of \(h\), Neumann proved that
both quantities are finite, and Pyber bounded both \(\zeta\)
and \(h\) exponentially. We determine \(h(n)\) for every positive
integer \(n\) and prove that \(h(n)=1+\sqrt{\zeta(n)}\) whenever
the constraint admits a nonabelian group. The covers are built
directly from the commutator geometry of a scalar quotient,
and the comparison with \(\zeta\) is made afterwards. Beyond
finitely many small cases, both maxima are attained by the same
groups, those isoclinic to an extraspecial \(2\)-group of order
\(2^{1+2\lfloor(n-1)/2\rfloor}\), confirming Pyber’s prediction
for the central index in this range.
\end{abstract}



\section{Introduction}\label{sec:intro}

How many abelian subgroups are needed to cover a group in which
every sufficiently large set contains two commuting elements?
The aim is to turn this restriction on noncommutation into
a sharp bound for the size of an abelian cover.
Write \(\omega(G)\) for the largest size of a set
of pairwise noncommuting elements, called a \emph{clique}, and
\(a(G)\) for the least size of an abelian cover. A third invariant,
the central index \(z(G)=[G:\Z(G)]\), measures the size of the
quotient on which commutation is determined. These satisfy
\[
 \omega(G)\le a(G)\le z(G):
\]
a clique meets each abelian subgroup in at most one point, and
each central coset \(x\Z(G)\) lies in \(\langle x,\Z(G)\rangle\).

Erd\H{o}s asked for the optimal bound on \(a(G)\) under the
constraint \(\omega(G)\le n\) \cite[Problem~26]{erdos}.
Neumann \cite{neumann} proved that this constraint bounds the
central index. Thus both extremal functions
\[
 h(n)=\max_{\omega(G)\le n}a(G),\qquad
 \zeta(n)=\max_{\omega(G)\le n}z(G)
\]
are finite. Pyber \cite{pyber} proved an exponential upper
bound for \(\zeta(n)\), hence also for \(h(n)\). Extraspecial
\(2\)-groups give the exponential lower bound attributed to
Isaacs \cite[p.~40]{bertram}.

The covering question also has a graph formulation. In the
noncommuting graph \(\Gamma_G\), the vertices are the noncentral
elements and edges join noncommuting pairs. Then
\(\omega(G)=\omega(\Gamma_G)\) for nonabelian \(G\), and
\(a(G)=\max\{1,\chi(\Gamma_G)\}\): a colour class generates an
abelian subgroup, and a cover gives a colouring.
Faber, Laver and McKenzie \cite{faber-laver-mckenzie} developed
this connection, while Mason \cite{mason} and Bertram
\cite{bertram} studied quantitative relations between covers
and noncommuting sets. See also Abdollahi, Akbari and Maimani
\cite{noncommuting-graph} for the graph-theoretic perspective.
Thus the covering problem asks for the optimal bound on the
chromatic number of noncommuting graphs in terms of their
clique number.

Here we determine \(h(n)\) for every \(n\), classify its
maximizers for \(n\ge9\), and identify its exact relation
to \(\zeta(n)\).
The covering profile and its equality case are proved directly.
The comparison with \(\zeta\) uses only the central-index theorem
\cite[Theorems~1.1 and~1.2(v)]{index}.
For \(n\ge9\), the maximizers all have binary symplectic
commutation geometry, in the dimension specified below.

\begin{theorem}[The covering profile]\label{thm:main}
For every \(n\ge3\),
\[
 h(n)=1+\sqrt{\zeta(n)}.
\]
Explicitly, \(h(n)=T(n)\) for all \(n\ge1\), where
\[
\begin{array}{c|cccc}
 n&1,2&3\le n\le6&7,8&n\ge9\\ \hline
 T(n)&1&n&10&1+2^{\lfloor(n-1)/2\rfloor}
\end{array}.
\]
For \(n\ge9\), the groups attaining \(h(n)\) under
\(\omega(G)\le n\) are exactly those isoclinic to an
extraspecial \(2\)-group of order \(2^{1+2\lfloor(n-1)/2\rfloor}\).
These are also the groups attaining \(\zeta(n)\).
\end{theorem}

\paragraph{Why a square root appears.}
The identity in Theorem~\ref{thm:main} reflects the geometry of
the extremizers. Let \(K=\F_{p^r}\). We call a group a
\emph{finite-field determinant model} if its central quotient
is the additive group of \(K^2\) and its commutator map is
\[
 \beta((x_1,x_2),(y_1,y_2))
 =\rho(x_1y_2-x_2y_1),
\]
where \(\rho\) is a nonzero \(\F_p\)-linear map onto the derived
subgroup, written additively. The determinant itself and its
scalar trace versions are both included. The \(K\)-lines give
a spread, and every such group satisfies
\[
 z(G)=|K|^2,\qquad a(G)=|K|+1=1+\sqrt{z(G)}
\]
(Corollary~\ref{cor:determinant-model}). The theorem shows that,
for each \(n\ge3\), a group of this type simultaneously
maximizes the covering number and the central index under
\(\omega(G)\le n\). In the small cases the relevant fields
have orders \(2,3,4,5,9\); for \(n\ge9\), the models are the
binary scalar trace models. Thus a common geometric construction
explains the relation between the two extremal functions.

For \(N\ge9\), put \(B(N)=2^{\lfloor(N-1)/2\rfloor}\), so
\(T(N)=1+B(N)\). The function \(B\) is nondecreasing, with \(B(N)\ge2^{N/2-1}\) and
\(B(N+2)=2B(N)\). In particular \(T(N)>N\): this holds at
\(N=9,10\) and persists in steps of two.

\paragraph{The idea of the proof.}
The proof compares cliques and covers in the fibres of a
quotient map. Call the largest clique size in a fibre its
\emph{width}. If the kernel is abelian, a fibre can be covered
by exactly that many abelian subgroups. More generally, an
abelian subgroup of index \(I\) in the kernel gives a cover
of size at most \(I\) times the width. A clique in the quotient
allows us to combine cliques in its fibres, so the sum of
their widths is at most \(\omega(G)\). Averaging over quotient
cliques then bounds the total number of subgroups in the cover
(Section~\ref{sec:clique-calculus}).

For a \(p\)-group, a scalar commutator form provides a symplectic
quotient. In rank two, an abelian radical already gives
\(a(P)=\omega(P)\) by Lemma~\ref{lem:rank-two-terminal}.
In the remaining cases, the symplectic group allows us to
average fibre cliques (Section~\ref{sec:primary}). When covering
the fibres separately gives too many subgroups, we instead use
subspaces isotropic for all commutator forms
(Section~\ref{sec:two}). The model is the symplectic spread,
which also constructs the sharp examples in
Section~\ref{sec:symplectic}. For nonnilpotent groups,
Section~\ref{sec:nonnil} uses affine geometry in a minimal
nonnilpotent quotient for the solvable case. A centralizer
descent applies to every nonnilpotent group when \(N\ge45\);
nonsolvable groups with \(N\le44\) are treated by the
small-clique classification.
For \(N\ge9\), each construction gives a strict inequality
except in the binary symplectic model with one commutator form.
Tracking where equality is possible therefore identifies
the maximizers.

The main distinction is between a scalar commutator form,
which certifies noncommutation, and the full commutator map,
which determines whether a subgroup is abelian. The former
provides large cliques; the latter provides common isotropic
subspaces and hence covers. Keeping both roles explicit allows
the same quotient estimates to handle groups of different
nilpotency classes. The case maps in
Sections~\ref{sec:primary}--\ref{sec:nonnil} indicate where each
estimate is used.

\paragraph{Classification input.}
The nonsolvable small-clique case uses the structural
classification in \cite[Theorem~1.1]{nonsol}; we do not claim
a proof independent of the classification of finite simple
groups (CFSG). The centralizer descent in
Section~\ref{sec:centralizer-descent} uses only the classical
Baer--Suzuki theorem.

\paragraph{Passing to finite groups.}
Two groups are \emph{isoclinic} if isomorphisms of their
central quotients and derived subgroups preserve the
commutators of lifts of central cosets.
A finite \(p\)-group is extraspecial when its centre and
derived subgroup coincide and have order \(p\).
For an extraspecial \(2\)-group the commutator is a
nondegenerate alternating form with values in its
derived subgroup of order two. Conversely, if \(G'\)
has order two and \(G/\Z(G)\) is elementary abelian
of order \(2^{2r}\), its commutator defines such a form.
An isometry to the form of an extraspecial group of
order \(2^{1+2r}\) gives an isoclinism. We call this the
\emph{binary symplectic model} of dimension \(2r\).

By Zarrin's theorem \cite[Theorem~1.1(2)]{zarrin}, a group of finite
clique number is isoclinic to a finite group. Isoclinism
preserves \(z\) and \(\omega\), since the central quotient
and the commutation of its lifts determine them. To see
that it preserves \(a\), enlarge each member of an abelian
cover by the centre. The resulting subgroups are exactly
the lifts of subgroups of the central quotient on which
the commutator is trivial. We may therefore work with
finite groups throughout.

Write \(N=\omega(G)\). A noncommuting pair \(x,y\)
extends to a clique \(x,y,xy\). Thus clique number two never
occurs, and clique number one characterizes abelian groups.
We use \([x,y]=x^{-1}y^{-1}xy\) and write \(C_m\) for a cyclic
group of order \(m\). The derived subgroup
is \(G'=[G,G]\); \(P^p\) is generated by the \(p\)-th powers.
Centralizers of subsets are taken pointwise. Vector-space
dimensions are over the indicated field; a line and a plane
have vector dimensions one and two. Projective dimensions are
specified explicitly. Orthogonality refers to the chosen scalar
form unless a subscript specifies another form. Logarithms
without a subscript are natural.

\section{Cliques and covers in cosets}\label{sec:clique-calculus}

For a subset \(X\) of a group, write \(\omega(X)\) for the
largest size of a clique contained in \(X\). We call
\(\omega(gU)\) the \emph{width} of the coset \(gU\).

\begin{lemma}[The abelian coset formula]\label{lem:abelian-coset}
If \(U\) is abelian, then
\begin{equation}\label{eq:abelian-coset}
 \omega(gU)=[U:C_U(g)].
\end{equation}
Moreover, \(gU\) can be covered by exactly that many
abelian subgroups.
\end{lemma}

\begin{proof}
For \(a,b\in U\), the elements \(ga,gb\) commute if and
only if \(g(ab^{-1})=(ab^{-1})g\), by cancelling the two
products and using commutativity in \(U\). Thus they
commute exactly when \(a,b\) lie in the same coset of
\(C_U(g)\) in \(U\). A clique takes at most one element
from each resulting class, and any such choice is a clique.
Each class consists of commuting elements and generates an
abelian subgroup. Conversely, the clique gives a lower bound
on the number of abelian subgroups in a cover.
\end{proof}

If \(H\) has an abelian subgroup of index \(I\), partition
\(gH\) into \(I\) cosets of that subgroup. Each has width
at most \(\omega(gH)\). By the lemma, each can be covered
by at most that many abelian subgroups. Thus \(gH\) costs at most
\(I\omega(gH)\).
Throughout, the \emph{cost} of a cover means its number of
abelian subgroups.

\subsection{Lifting a subgroup cover}

For nonabelian \(G\), fix an abelian cover of a subgroup
\(U\le G\). For each coset \(gU\ne U\), its translates
\(gA\), with \(A\) a member of this cover, cover \(gU\).
Each translate costs at most \(N-1\) subgroups. If \(g\)
fails to centralize \(A\), an element of \(A\) extends a
maximum clique in \(gA\). If \(g\) centralizes \(A\), the
translate costs one. Thus
\begin{equation}\label{eq:lift}
 a(G)\le a(U)\bigl[1+([G:U]-1)(N-1)\bigr].
\end{equation}

\subsection{Adding widths across fibres}

One fibre per nontrivial cyclic subgroup of the quotient
suffices for the cover. To control their total width, we
average over families whose cliques can be joined.

\begin{lemma}[Clique and cover bounds in a quotient]\label{lem:fibre-calculus}
Let \(\pi:G\to Q\) be a surjection with kernel \(H\).
Write \(w(q)=\omega(\pi^{-1}(q))\), and let \(I\) be the
index of an abelian subgroup of \(H\).
\begin{enumerate}
\item[(i)] Two generators of the same cyclic subgroup of
\(Q\) have fibres of the same width, covered by the same
families of abelian subgroups.
\item[(ii)] If every point of \(X_i\subseteq Q\) fails to
commute with every point of \(X_j\subseteq Q\) for \(i\ne j\), then
\[
 \sum_i\omega(\pi^{-1}(X_i))\le\omega(G).
\]
\item[(iii)] If \(D\subseteq Q\) contains a generator of every
nontrivial cyclic subgroup of \(Q\), then
\[
 a(G)\le a(H)+I\sum_{q\in D}w(q).
\]
Suppose a random subset of \(D\) contains each \(q\in D\)
with probability at least \(\delta>0\). If lifts of its
distinct points are always noncommuting, then
\(\sum_{q\in D}w(q)\le\omega(G)/\delta\).
\end{enumerate}
\end{lemma}

\begin{proof}
For (i), choose an exponent inducing the required power on
the cyclic subgroup and coprime to \(|G|\), using the Chinese
remainder theorem. The power map is a bijection between the
fibres and preserves the cyclic subgroup generated by each
element. This preserves both commutation and membership in
abelian subgroups. For (ii), take a maximum clique in
each inverse image. Their union is a clique in \(G\).
For (iii), cover \(H\) and each selected fibre, then use
(i). For every realization of the random subset, the
cliques in its fibres can again be joined. Their total
size is at most \(\omega(G)\), and its expectation is at
least \(\delta\sum_{q\in D}w(q)\).
\end{proof}

Lifting the cyclic subgroups of \(G/\Z(G)\) gives
\(a(G)\le\max(1,z(G)-1)\). For a \(p\)-group of
central index \(p^k>1\) this cost is at most
\[
\mathcal C_p(k):=\frac{p^k-1}{p-1},
\]
since each nontrivial cyclic subgroup has at least \(p-1\) generators.

\subsection{Centralizers in a narrow coset}

The factor \(I\) makes the quotient estimate useful only
if the kernel has a large abelian subgroup. To find one, we
look inside a narrow coset. The centralizers of a maximum
clique cover that coset, so one of them meets a large part
of the kernel. The same cover also controls the intersection
of the centralizers.

\begin{lemma}[Centralizer partitions]\label{lem:centralizer-partition}
Let \(D\) be a maximum clique in a finite subset \(X\) of \(G\).
The sets \(X\cap C_G(x)\), \(x\in D\), cover \(X\).
If the sum of their sizes is at most \(|X|\), they partition
\(X\) into commuting sets.
\end{lemma}

\begin{proof}
An uncovered point would extend \(D\). If the sizes sum to
at most \(|X|\), the cover is a partition. A noncommuting
pair in the part at \(x\) would replace \(x\) and enlarge
\(D\), since neither point commutes with any other vertex.
\end{proof}

We call a cover \emph{irredundant} if deleting any member
leaves an uncovered point.

\begin{lemma}[Centralizers in a coset]\label{lem:coset-calculus}
Let \(U\le G\), \(g\in G\), and let \(X\) be a maximum clique in \(gU\).
Then
\[
 U\cap C_G(X)=\Z(U)\cap C_U(g).
\]
The sets \(C_G(x)\cap gU\), \(x\in X\), form an irredundant coset
cover. Each \(x\in gU\) satisfies \(\omega(C_U(x))\le\omega(gU)\),
and some \(x\in X\) satisfies \([U:C_U(x)]\le|X|\).
If \(U\) is a nonabelian \(p\)-group, then \(\omega(gU)\ge p+1\).
\end{lemma}

\begin{proof}
Suppose \(a\in U\) centralizes \(X\) but fails to
commute with \(y\in gU\). Replace each vertex \(x\)
that commutes with \(y\) by \(xa\). The replacements
remain in the coset and preserve mutual commutators
because \(a\) centralizes \(X\). Now no vertex commutes
with \(y\), so we can adjoin \(y\), a contradiction.
Thus \(a\) centralizes all of \(gU\), which is equivalent
to belonging to \(\Z(U)\cap C_U(g)\).

The parts of the cover in Lemma~\ref{lem:centralizer-partition}
are \(xC_U(x)\), and each vertex belongs only to its own part.
At least one part has size \(|U|/|X|\) or more, giving the
index bound. Multiplication by \(x\) carries a clique in
\(C_U(x)\) into \(gU\) with its commutation relations
unchanged. Translating the cover by \(g^{-1}\) on the left
and then inverting gives an irredundant cover of \(U\)
by right cosets of the \(C_U(x)\).

If \(U\) is a nonabelian \(p\)-group, a one-part cover
would contradict the intersection formula. Each subgroup
in the cover is therefore proper. At most \(p\) parts
could cover only by partitioning \(U\) into \(p\)
cosets of index \(p\). Two distinct index-\(p\)
subgroups are normal and multiply to \(U\), so any two
of their cosets intersect. Hence all the subgroups in
the partition would coincide. Their intersection lies
in \(\Z(U)\), which would make \(U/\Z(U)\) cyclic.
\end{proof}

Tomkinson's bound for irredundant coset covers
\cite{tomkinson} now gives
\begin{equation}\label{eq:tomkinson-coset}
 [U:\Z(U)\cap C_U(g)]\le\omega(gU)! .
\end{equation}
In particular a coset of a nonabelian subgroup has width
at least three.

Suppose \(x\in gU\) centralizes an abelian subgroup \(A\)
of index \(I\) in \(U\). We can use the fibre width to cover
\(U\) itself. For each \(u\in U\), the elements
of \(A\) commuting with \(u\) are exactly those commuting
with \(xu\). The abelian coset formula then gives
\[
 \omega(uA)=[A:C_A(u)]=[A:C_A(xu)]=\omega(xuA)\le\omega(gU).
\]
Covering the \(I\) cosets of \(A\) in \(U\) therefore costs
at most \(I\omega(gU)\). 

\section{Symplectic geometry and the sharp examples}\label{sec:symplectic}

When an alternating form determines commutation in a central
quotient, cliques bound the admissible profile and isotropic
subspaces lift to abelian subgroups. We construct both.

An alternating form \(b\) is \emph{symplectic} if the
only vector pairing to zero with every vector is zero.
A clique for \(b\) is a set whose distinct vectors pair
nontrivially. In dimension two, the basic example is
the determinant pairing. Two nonzero vectors pair to
zero exactly when they lie on the same line. Thus a
clique is obtained by choosing one nonzero vector on
each line. The following lemma extends this count to
larger dimensions. For odd primes, the lower bound \(pr+1\)
is due to Chin \cite{chin}, in the equivalent formulation
for extraspecial \(p\)-groups. The binary maximum is due to
Isaacs \cite[p.~40]{bertram}. We include proofs for completeness.

\begin{lemma}[Symplectic cliques]\label{lem:symplectic-cliques}
Let \((V,b)\) be symplectic of dimension \(2r>0\) over \(\F_p\).
There is a clique of size \(pr+1\).
For \(p=2\), the maximum is \(2r+1\), and every maximum
clique spans \(V\) and has sum zero.
\end{lemma}

\begin{proof}
To combine cliques \(X,Y\) in orthogonal summands,
choose \(y_0\in Y\) and use
\((X+y_0)\cup(Y\setminus\{y_0\})\).
The translate retains its original pairings, and its
pairings with the remaining vertices of \(Y\) are those
of \(y_0\). The new clique has size \(|X|+|Y|-1\).
Starting with \(p+1\) vertices per symplectic plane gives
\(pr+1\) vertices.

Over \(\F_2\), a clique of size \(k\) has a Gram matrix
with diagonal zero and all other entries one. Its
kernel consists of the vectors whose entries all
equal their sum. The rank is therefore \(k\) for even
\(k\) and \(k-1\) for odd \(k\). It cannot exceed \(2r\),
so \(k\le2r+1\). At equality the rank is \(2r\), so the
clique spans \(V\). Its sum is orthogonal to the clique
and is therefore zero.
\end{proof}

We call a maximum clique for a binary symplectic form
a \emph{symplectic simplex}.

A subspace on which an alternating form vanishes on
every pair is called \emph{totally isotropic}. When
commutation in a central quotient is determined by
such a form, the abelian subgroups containing the
centre are exactly the lifts of these subspaces.

\begin{lemma}[A cover from a symplectic spread]\label{lem:spread}
Suppose the central quotient of \(G\) is a vector space
of dimension \(2r>0\) over \(\F_p\), and commutation
of lifts is determined by a symplectic form on that
space. Then \(a(G)=p^r+1\).
\end{lemma}

\begin{proof}
A totally isotropic subspace \(U\) satisfies \(U\le U^\perp\),
so \(\dim U\le r\). It contains at most \(p^r-1\) nonzero
vectors. Dividing the total \(p^{2r}-1\) by \(p^r-1\)
gives the lower bound \(p^r+1\).

For the cover, write \(K=\F_{p^r}\) and view \(K^2\) as
a space over \(\F_p\). Use the form
\[
 ((x_1,x_2),(y_1,y_2))\longmapsto
 \Tr_{K/\F_p}(x_1y_2-x_2y_1).
\]
The trace pairing is nondegenerate, so this form is symplectic.
Every \(K\)-line is totally isotropic because the
determinant vanishes on it. The \(p^r+1\) lines
partition the nonzero vectors of \(K^2\). Such a
partition into totally isotropic subspaces of
dimension \(r\) is called a \emph{symplectic spread}.
All nondegenerate alternating forms over \(\F_p\)
of dimension \(2r\) have symplectic bases, so the
partition gives a spread in the original quotient.
Lifting its members gives the required cover.
\end{proof}

\begin{corollary}[Finite-field determinant models]
\label{cor:determinant-model}
Every finite-field determinant model with central quotient
\(K^2\) satisfies \(a(G)=|K|+1=1+\sqrt{z(G)}\).
\end{corollary}

\begin{proof}
Write \(|K|=p^r\). The \(|K|+1\) lines over \(K\) are
isotropic for the determinant and hence for its image under
\(\rho\). Their preimages form an abelian cover.
For the converse, choose a linear functional \(\lambda\) on
the derived subgroup such that \(\lambda\rho\ne0\). The form
\(\lambda\rho\circ\det\) is nondegenerate over \(\F_p\): for
\(v\ne0\), the map \(w\mapsto\det(v,w)\) is onto \(K\).
Every abelian subgroup, enlarged by the centre, therefore
has isotropic image of dimension at most \(r\). Counting
nonzero vectors as in Lemma~\ref{lem:spread} requires at least
\((p^{2r}-1)/(p^r-1)=p^r+1\) members.
\end{proof}

\subsection{Groups attaining the bounds}\label{sec:models}

The constructions above give all lower bounds in
Theorem~\ref{thm:main}. Let \(k\) be the allowed clique size.

For \(k=1,2\) the abelian model attains the bound.
For \(3\le k\le6\), take \(\mathrm{UT}_3(\F_q)\), the group of
upper unitriangular matrices, with \(q=k-1\in\{2,3,4,5\}\).
Its central quotient is \(\F_q^2\), and the commutator is the
determinant \(xy'-x'y\). The \(q+1\) lines lift to an abelian
cover; one nonzero point on each line gives a clique of size
\(q+1=k\).

For \(k=7,8\), use an extraspecial \(3\)-group with a
four-dimensional symplectic central quotient. By the duality
between the symplectic generalized quadrangle \(W(3,q)\) and
the parabolic quadric \(Q(4,q)\), a set of pairwise nonorthogonal
projective points corresponds to a partial spread of the quadric.
Tallini's bound \cite{tallini}\footnote{Klein and Metsch cite
Tallini's 1991 paper \cite{tallini} for the odd-\(q\) result;
Ceria et al.~\cite[Introduction]{ceria} attribute the bound to
his earlier paper \cite{tallini-fibrations}. We use the
partial-spread formulation in Klein--Metsch.}, together with the nonexistence
of ovoids in \(W(3,q)\) for odd \(q\), bounds its size by
\(q^2-q+1\); see \cite[Theorem~3.3 and Section~6]{klein-metsch}.
At \(q=3\) this is seven, attained by combining two plane
cliques. The clique number is therefore seven, and the spread
in Lemma~\ref{lem:spread} gives an optimal cover of size ten.

For \(k\ge9\), take an extraspecial \(2\)-group of order
\(2^{1+2r}\) with \(r=\lfloor(k-1)/2\rfloor\). Its central
index is \(4^r\), its clique number is \(2r+1\), and its
covering number is \(2^r+1\), by
Lemmas~\ref{lem:symplectic-cliques} and~\ref{lem:spread}.

\section{The scalar quotient of a \texorpdfstring{\(p\)}{p}-group}\label{sec:primary}

A finite nilpotent group is a product of its Sylow subgroups.
We first cover each prime-power factor using a scalar
commutator form. We will combine their covers at the end of
Section~\ref{sec:two}.

\paragraph{Case map.}
For a scalar radical \(R\), write \(d=\dim(P/R)=2j\).
The binary cases are treated in Section~\ref{sec:two}; the
odd-prime cases below are organized as follows. In the last
row, \([R:\Z(R)]=p^m\).
\begin{center}
\renewcommand{\arraystretch}{1.12}
\begin{tabular}{@{}p{.25\textwidth}p{.23\textwidth}p{\dimexpr0.52\textwidth-4\tabcolsep\relax}@{}}
\hline
\casecell{Radical and prime} & \casecell{Range or constraint} & \casecell{Covering argument} \\
\hline
\casecell{Any radical, \(p=2\)} & \casecell{All \(N\)} & \casecell{Binary case map in Section~\ref{sec:two}.} \\
\casecell{Abelian noncentral, \(p\) odd} & \casecell{\(j=1\); or \(N\ge p^2j+1\)} & \casecell{Lemma~\ref{lem:rank-two-terminal}; otherwise the fibre average in Section~\ref{sec:fibre-average}.} \\
\casecell{Central, one component, \(p\) odd} & \casecell{\(j=1\); or \(p\ge5\), \(j\ge2\)} & \casecell{The spread of Lemma~\ref{lem:spread}, with \(N\ge pj+1\).} \\
\casecell{Central, one component, \(p=3\)} & \casecell{\(j=2\); or \(j\ge3\)} & \casecell{The ten-member cover; otherwise thirteen-point cliques and the spread.} \\
\casecell{Central, several components, \(p\) odd} & \casecell{\(N\ge p^2(j-1)+1\)} & \casecell{Lemma~\ref{lem:affine-inflation}; the common centralizer planes in dimension four.} \\
\casecell{Nonabelian, \(p\) odd} & \casecell{\(N\ge(pj+1)\max\{p+1,m+1\}\)} & \casecell{Lemma~\ref{lem:radical-width} and the fibre average.} \\
\hline
\end{tabular}
\end{center}

\begin{proposition}[Nilpotent groups]\label{prop:nilpotent-cover}
If \(G\) is nilpotent and \(N=\omega(G)\), then \(a(G)\le T(N)\).
For \(N\ge9\), equality requires the binary symplectic
model of dimension \(N-1\).
\end{proposition}

\subsection{The quotient given by one commutator form}
\label{sec:scalar-fibres}

Let \(P\) be a nonabelian \(p\)-group and set \(N=\omega(P)\).
A single scalar commutator value will suffice to certify
noncommutation of lifts.

The commutators are central modulo
\(\gamma_3(P)=[P',P]\). The abelian \(p\)-group
\(P'/\gamma_3(P)\) is nontrivial: otherwise the lower
central series would be constant from \(P'\), contrary
to nilpotence. Choose a nonzero homomorphism
\(\xi:P'/\gamma_3(P)\to(\F_p,+)\), and also write \(\xi\)
for its composite with the quotient map on \(P'\).

Applying \(\xi\) to commutators gives the pairing
\[
 b_\xi(x,y)=\xi([x,y]).
\]
Centrality of the commutators in \(P/\gamma_3(P)\)
makes the pairing additive in each argument. Its
\emph{radical} is the subgroup
\[
 R=R_\xi=\{x\in P:b_\xi(x,y)=0\text{ for all }y\in P\}.
\]
It is normal and contains \(P'P^p\Z(P)\). Hence
\(V=P/R\) is a vector space over \(\F_p\), and the
induced form is symplectic by the definition of \(R\).
We call \(b_\xi\) a \emph{scalar commutator form} and
\(R\) its \emph{scalar radical}. A nonzero pairing in \(V\) certifies noncommutation of every
pair of lifts, so a scalar clique lets us add its fibre widths.

Scalar radicals need not be central, even in class two:
\(R_\xi\) tests only \(\xi([x,y])=0\), whereas \(\Z(P)\)
requires \([x,y]=1\) for every \(y\).
For example, in \(D_8\times D_8\), projecting the commutator
onto the first component gives the nonabelian radical
\(\Z(D_8)\times D_8\); here \(D_8\) is the dihedral group
of order eight.

\begin{lemma}[The terminal rank-two case]\label{lem:rank-two-terminal}
If a scalar commutator form of \(P\) has rank two and an
abelian radical \(R\), then \(a(P)=\omega(P)\).
\end{lemma}

\begin{proof}
For \(x\notin R\), the perpendicular of \(xR\) in the symplectic
plane \(P/R\) is its own line. Thus
\(C_P(x)=\langle x,C_R(x)\rangle\) is abelian.
Two commuting elements outside \(R\) have the same centralizer:
each centralizer is abelian and hence contained in the other.

Choose one representative outside \(R\) for each distinct centralizer.
They form a clique, and their centralizers cover \(P\setminus R\).
If they cover \(R\) as well, the clique and cover have the same size.
Otherwise an uncovered point of \(R\) extends the clique,
and adding \(R\) extends the cover. Hence \(a(P)=\omega(P)\).
\end{proof}

\subsection{Adding the widths of the fibres}
\label{sec:fibre-average}

Write \(d=2j\) for the dimension of \(V\). Select one nonzero
vector on each line \(L\) of \(V\).
Its fibre has width \(w_L\), and a cover of that fibre
serves the whole line by Lemma~\ref{lem:fibre-calculus}(i).
To bound the total width, take the lines of a symplectic
clique of size \(pj+1\) and move them by a uniformly chosen
symplectic transformation. The selected vectors on the
resulting lines still form a clique. Transitivity on lines
gives inclusion probability \((pj+1)/\mathcal C_p(2j)\).
Lemma~\ref{lem:fibre-calculus}(iii) now gives
\[
 \sum_L w_L\le\frac{\mathcal C_p(2j)\,N}{pj+1}.
\]
If \(R\) has an abelian subgroup of index \(I\), the same
lemma gives the common bound
\begin{equation}\label{eq:uniform-fibre-cost}
 a(P)\le a(R)+I\,\frac{\mathcal C_p(2j)\,N}{pj+1}.
\end{equation}

\begin{lemma}[A clique meeting a hyperplane once]\label{lem:clique-meets-hyperplane}
In a symplectic space \((V,b)\) of dimension \(2r>0\) over \(\F_p\),
every hyperplane meets some clique of size \(pr+1\) in
exactly one point.
\end{lemma}

\begin{proof}
For a prescribed hyperplane \(e^\perp\), choose \(f\)
with \(b(e,f)=1\), and let \(V_0\) be the orthogonal
complement of \(\langle e,f\rangle\). Take a clique
\(D_0\) of size \(p(r-1)+1\) in \(V_0\), with
\(D_0=\{0\}\) if \(r=1\). Then
\[
 (D_0+f)\cup\{e+tf:t\in\F_p\}
\]
is a clique. The first part inherits the pairings of
\(D_0\). In the second part they are differences of
distinct coefficients of \(f\), and the cross pairings
are all \(b(f,e)\). Only \(e\) lies in \(e^\perp\).
\end{proof}

\begin{lemma}[Coset covers of \(p\)-groups]\label{lem:coset-cover}
If \(U\) is a finite \(p\)-group and the right cosets
\(H_i x_i\), \(1\le i\le t\), form an irredundant cover
of \(U\), then \([U:\bigcap_i H_i]\le p^{t-1}\).
\end{lemma}

\begin{proof}
We argue by induction on \(t\). At each step, we restrict the
cover to a coset disjoint from the members already selected.
The restricted cover uses fewer members, so induction controls
the index of the next intersection. The case \(t=1\) is immediate.

For \(t>1\), all \(H_i\) are proper. Choose a maximal subgroup
\(M\) containing \(H_1\). It is normal of index \(p\).
Right-translate the cover so that \(H_1x_1\) misses \(M\),
and choose an irredundant subcover of \(M\), with index set
\(I_1\). Put \(B_1=\bigcap_{i\in I_1}H_i\). Since \(1\notin I_1\),
induction gives \([M:M\cap B_1]\le p^{|I_1|-1}\).

Write \(D=\bigcap_iH_i\). If the intersection \(B_j\) of the
subgroups selected so far is not yet \(D\), the selected
cosets do not cover \(U\), by irredundancy. Choose a point \(y\)
that they miss. They miss the whole coset \(B_jy\), because
\(B_j\) lies in every selected subgroup. Restrict the remaining
cover to \(B_jy\) and translate by \(y^{-1}\). An irredundant
subcover, indexed by \(I_{j+1}\), consists of cosets of
\(B_j\cap H_i\) in \(B_j\). Its indices are all new, so induction
bounds the next intersection by
\[
 [B_j:B_{j+1}]\le p^{|I_{j+1}|-1},\qquad
 B_{j+1}=B_j\cap\bigcap_{i\in I_{j+1}}H_i.
\]
Each step adds a nonempty set of previously unused indices,
so the process terminates at \(B_s=D\). Since \(D\le M\),
multiplying the index bounds gives
\[
 [U:D]\le p[M:M\cap B_1][B_1:D]
 \le p^{1+\sum_{j=1}^s(|I_j|-1)}.
\]
If the index \(1\) remains unused, then
\(\sum_j|I_j|\le t-1\) and \(s\ge1\), so the exponent
is at most \(t-1\). If it is used, then \(s\ge2\), since
\(1\notin I_1\), and \(\sum_j|I_j|\le t\) gives the same bound.
\end{proof}

\begin{lemma}[Fibres over a noncentral radical]\label{lem:radical-width}
Let \(R\) be the radical of a scalar commutator form of rank
\(2j\) on \(P\), and put \(\sigma=pj+1\).
\begin{enumerate}
\item If \(R\) is abelian and noncentral, then \(\omega(P)\ge p^2j+1\).
\item If \(R\) is nonabelian and \([R:\Z(R)]=p^m\), then the width
of each fibre \(gR\) is at least \(p+1\) and at least
\(m+1+\log_p[\Z(R):C_{\Z(R)}(g)]\). Consequently
\(\omega(P)\ge \sigma\max\{p+1,m+1\}\), and
\(\omega(P)\ge \sigma(m+1)+1\) if \(\Z(R)\ne\Z(P)\).
\end{enumerate}
\end{lemma}

\begin{proof}
For (1), choose \(a\in R\setminus\Z(P)\) and a hyperplane of \(P/R\)
containing \(C_P(a)/R\). By Lemma~\ref{lem:clique-meets-hyperplane},
some symplectic clique of size \(\sigma\) meets the hyperplane once.
Each of its \(pj\) other vertices has a fibre of width at least \(p\)
by \eqref{eq:abelian-coset}, since every lift moves \(a\). The remaining
fibre contributes at least one.

For (2), let \(w\) be the width of \(gR\). Lemma~\ref{lem:coset-calculus}
gives \(w\ge p+1\) and an irredundant cover of \(R\) by \(w\) cosets.
Their subgroups intersect in \(\Z(R)\cap C_R(g)\), of index at most
\(p^{w-1}\) by Lemma~\ref{lem:coset-cover}. This index is also
\(p^m[\Z(R):C_{\Z(R)}(g)]\), giving the second width bound. Summing
the widths over a clique of size \(\sigma\) gives
\(\omega(P)\ge \sigma\max\{p+1,m+1\}\).
If \(\Z(R)\ne\Z(P)\), choose \(g\) moving a point of \(\Z(R)\).
Then \(g\notin R\) and \(gR\) has width at least \(m+2\).
Some symplectic clique of size \(\sigma\) contains \(gR\), by
transitivity on nonzero vectors, giving the last bound.
\end{proof}

\subsection{Odd primes}

Choose a scalar form of least nonzero rank. The fibre average
will suffice except in the central ternary cases, where we use
larger cliques or a common isotropic spread.

\paragraph{Abelian noncentral radical.}
Lemma~\ref{lem:rank-two-terminal} gives \(a(P)=N\) when \(j=1\).
For \(j\ge2\), Lemma~\ref{lem:radical-width} gives \(N\ge p^2j+1\).
To compare the fibre-cover bound with the profile, note that
\(N2^{-N/2}\) decreases for \(N\ge3\),
since each step multiplies it by \((1+1/N)/\sqrt2<1\).
Using \(B(N)\ge2^{N/2-1}\), we may compare at \(N=p^2j+1\).
There \((p^2j+1)/(pj+1)<p\), so
\[
 \frac{a(P)-1}{B(N)}
 < \frac{\sqrt2\,p}{p-1}
   \left(\frac{p^2}{2^{p^2/2}}\right)^j<1.
\]
The factor in parentheses is at most \(1/2\).
For \(p=3\) this follows from \(9<2^{7/2}\).
For \(p\ge5\) use \(p^2<2^p\le2^{p^2/2-1}\).
The first inequality holds at \(5\) and persists in steps
of two, since \(p+2<2p\).

\paragraph{Central radical and one commutator component.}
If \(R=\Z(P)\), then \(P\) has class two and \(P/\Z(P)\)
is elementary abelian. Its commutator subgroup is elementary
abelian as well, because \([x,y]^p=[x^p,y]=1\). Write
\(\beta\) for the alternating commutator map on \(P/\Z(P)\)
with values in \(P'\). Composing \(\beta\) with a linear
functional on \(P'\) gives a \emph{scalarization} of \(\beta\).
When \(P'\) has dimension one, a symplectic
spread gives exactly \(p^j+1\) subgroups.
In a scalar plane the spread costs exactly \(N=p+1\).
For \(j\ge2\) and \(p\ge5\), we compare the spread size
with the profile at the scalar clique bound.
If \(j\) is even, \(p^2<2^p\) gives \(p^j<B(pj+1)\).
If \(j\) is odd, \(j\ge3\) and
\[
 \frac{p^j}{B(pj+1)}
 =\sqrt{2}\,(p\,2^{-p/2})^j
 \le \sqrt{2}\,(5\,2^{-5/2})^3=\frac{125}{128}.
\]
Here \(p\,2^{-p/2}\) decreases as \(p\) increases in steps
of two.

At \(p=3\) and \(j=2\), the spread costs ten and the clique
number is seven by Section~\ref{sec:models}. For \(j\ge3\),
we need larger cliques than the plane count provides.
Ceria, De Beule, Pavese and Smaldore \cite[Theorem~3.7]{ceria}
give thirteen projective points forming a clique in
dimension six over \(\F_3\): more generally, for
\(c\in\F_q^\times\), their \(q^2+q+1\) points
\[
 P_x=(x,x^q,x^{q^2},cx^{-q^2},cx^{-q},cx^{-1}),
 \qquad x\in\F_{q^3}^\times,\quad
 \operatorname{Nm}_{\F_{q^3}/\F_q}(x)=1,
\]
lie in the \(\F_q\)-subgeometry
\(\{[a:a^q:a^{q^2}:b^{q^2}:b^q:b]:(a,b)\ne(0,0)\}\),
where \(a,b\in\F_{q^3}\),
of \(\mathrm{PG}(5,q^3)\) and are pairwise nonorthogonal for
the alternating form
\(b(v,w)=\sum_{i=1}^3(v_iw_{7-i}-v_{7-i}w_i)\), since
\[
 b(P_x,P_y)
 =c\Tr_{\F_{q^3}/\F_q}(x/y-y/x)
 =c\operatorname{Nm}_{\F_{q^3}/\F_q}(1-y/x)\ne0
 \qquad(x\ne y).
\]
Taking \(q=3\) gives the required thirteen points.
Write \(j=3u+v\), with
\(0\le v\le2\), and combine \(u\) such cliques with \(v\)
plane cliques by Lemma~\ref{lem:symplectic-cliques}.
This gives \(N\ge12u+3v+1\ge4j-1\). Thus \(B(N)\ge4^j/2>3^j\),
since \(4^3/2>3^3\) and the bound grows faster than the
spread cost.

\paragraph{Central radical and several commutator components.}
A subspace isotropic for one form may contain points with
noncommuting lifts. The cover must therefore use subspaces
isotropic for every commutator form. To compensate for any
increase in their number,
we use a second scalarization to construct a larger clique. The
first form separates the fibres, and the second lets us
take a whole line of points in most of them.

\begin{lemma}[Cliques from two commutator forms]
\label{lem:affine-inflation}
Let \(\beta:V\times V\to W\) be alternating bilinear over \(\F_p\),
with values spanning \(W\), \(\dim V=2r\) and \(\dim W\ge2\).
If some scalarization \(b=\lambda\circ\beta\) is nondegenerate,
then \(r\ge2\) and there are \(p^2(r-1)+1\) vectors whose
\(\beta\)-values on distinct pairs are all nonzero.
\end{lemma}

\begin{proof}
Take an independent functional \(\mu\) on \(W\) and
put \(c=\mu\circ\beta\). Write \(c(v,x)=b(Sv,x)\).
The operator \(S\) is nonscalar because the values of
\(\beta\) span \(W\). Thus some \(v,Sv\) are independent, since an operator preserving
every line is scalar.

On the symplectic quotient \(v^{\perp_b}/\langle v\rangle\),
the functional \(c(v,\cdot)\) is well defined and
nonzero, since \(Sv\notin\langle v\rangle\). This also
shows that the quotient has positive dimension, so
\(r\ge2\). Choose a clique of size \(p(r-1)+1\)
meeting its kernel once, by
Lemma~\ref{lem:clique-meets-hyperplane}. Lift its vertices
to \(v^{\perp_b}\) and replace every lift \(x\) outside
the kernel by the line \(x+\F_pv\). The \(b\)-pairings
between different fibres are still nonzero. Within a
line the \(c\)-pairing is a nonzero multiple of
\(c(v,x)\), so all its points can be used.
This gives \(1+p^2(r-1)\) vertices.
\end{proof}

Applying the lemma to \(\beta\) gives \(N\ge p^2(j-1)+1\).
Lifting the lines gives a cover of size \(\mathcal C_p(2j)\).
For \(p\ge5\) and \(j\ge3\), the inequality \(p^2<2^p\) gives
\(p^{2j}<2^{pj}\). Since \(j\ge2\) and \(p\ge5\),
\(p^2(j-1)\ge2pj\), so the line count is below \(B(N)\).
For \(p=3,j\ge3\), the bound \(N\ge9(j-1)+1\) gives
\(B(N)\ge 2^9\,16^{j-3}\), while the line count is
less than \((9^3/2)9^{j-3}\).

Let \(j=2\). The least-rank choice implies that every nonzero
scalarization of \(\beta\) on the four-space is nondegenerate.
The Pfaffian bound \cite[Lemma~3]{gow} gives \(\dim P'\le2\),
and hence \(\dim P'=2\).
For \(x\ne0\), the map \(\beta(x,\cdot):V\to P'\) is onto:
a functional annihilating its image would give a nonzero
scalarization with \(x\) in its radical. Thus \(C_\beta(x)=
\ker\beta(x,\cdot)\) is a plane. It is isotropic for \(\beta\),
since it is spanned by \(x\) and a vector commuting with \(x\).
For each nonzero \(y\in C_\beta(x)\) we then have
\(C_\beta(x)\le C_\beta(y)\), with equality by dimension.
The planes partition the nonzero vectors of \(V\); their
\(p^2+1\) lifts give an abelian cover. Since \(N\ge p^2+1\ge10\)
by Lemma~\ref{lem:affine-inflation}, this cover is strictly below
the profile.

\paragraph{Nonabelian radical.}
Write \([R:\Z(R)]=p^m\) with \(m\ge2\), and put
\(\sigma=pj+1\) and \(M=\max\{p+1,m+1\}\).
Lemma~\ref{lem:radical-width} gives \(N\ge\sigma M\).
At most \(\mathcal C_p(m)\) cyclic lifts cover \(R\), and
one such lift has index at most \(p^{m-1}\) in \(R\).
Since \(\mathcal C_p(m)<\tfrac32p^{m-1}\le p^{m-1}N/\sigma\),
\eqref{eq:uniform-fibre-cost} gives
\[
 a(P)\le \mathcal C_p(m)+p^{m-1}\frac{\mathcal C_p(2j)N}{\sigma}
 \le p^{m-1}\bigl(1+\mathcal C_p(2j)\bigr)\frac{N}{\sigma}.
\]
The ratio of the last bound to \(2^{N/2-1}\) is largest at
\(N=\sigma M\), because \(N2^{-N/2}\) decreases. Since
\(1+\mathcal C_p(2j)\le(p+2)p^{2j-2}\), this gives
\[
 \frac{a(P)}{B(N)}\le \rho\,\theta^{j-1},\qquad
 \rho=2M(p+2)p^{m-1}2^{-(p+1)M/2},\qquad
 \theta=p^2\,2^{-pM/2}.
\]
Here \(\theta\le9/64\): \(M\ge4\), and \(p^2/4^p\) decreases
for odd \(p\ge3\).
For \(p\ge5\), use \(m\le M-1\) and
\(p\le\tfrac58\,2^{(p+1)/2}\), which holds at \(5\) and
persists in steps of two. Then
\[
 \rho\le\frac{2(p+2)}{p^2}M(5/8)^M
 \le\frac{14}{25}\,6(5/8)^6<1,
\]
since the two factors decrease for \(p\ge5\) and \(M\ge6\),
respectively. If \(p=3,m=2\), we have \(\rho=120/2^8<1\).
If \(p=3,m\ge3\), then \(M=m+1\) and
\(\rho=(10/9)M(3/4)^M\le45/32\), since \(M(3/4)^M\)
decreases for \(M\ge4\). This suffices when \(j\ge2\),
as \((45/32)(9/64)<1\).

It remains to take \(p=3,j=1,m\ge3\). Here
\(\sigma=\mathcal C_3(2)=4\), and the first cost bound gives
\(a(P)<3^{m-1}(N+2)\).
If \(\Z(R)\ne\Z(P)\), Lemma~\ref{lem:radical-width} gives
\(N\ge4m+5\). At \(N=4m+6\), the inequality
\(3^{m-1}(N+2)\le B(N)\) reduces to \(3^{m-1}(m+2)\le4^m\).
This holds at \(m=3\) and persists since \(3(m+3)\le4(m+2)\).
It also holds at \(N=4m+5\), where \(B\) has the same value,
and persists in steps of two, since \(B\) doubles.
If \(\Z(R)=\Z(P)\), then \(z(P)=3^{m+2}\) and cyclic
lifts give
\[
 a(P)\le\mathcal C_3(m+2)<3^{m+2}/2\le 2^{2m+1}\le B(N).
\]
The middle inequality starts at \(m=3\) and persists as
\(m\) increases; the last uses \(N\ge4m+4\).

Thus \(a(P)\le T(N)\) for odd \(p\). When \(N\ge9\), the
bounds above are all strict, including the costs \(N\)
and ten.

\section{Binary groups and the equality model}\label{sec:two}

\paragraph{Case map.}
Write \(\nu=\log_2 z(P)\). First use cyclic lifts when
\(N\ge2\nu+1\). In the remaining range \(N\le2\nu\),
test the radicals in the order shown below. Here \(d=\dim(P/R)\).
\begin{center}
\renewcommand{\arraystretch}{1.12}
\begin{tabular}{@{}p{.26\textwidth}p{.24\textwidth}p{\dimexpr0.50\textwidth-4\tabcolsep\relax}@{}}
\hline
\casecell{Radical or commutator} & \casecell{Range} & \casecell{Covering argument} \\
\hline
\casecell{Any} & \casecell{\(N\ge2\nu+1\)} & \casecell{Cyclic lifts; the small values \(N\le8\) are included below.} \\
\casecell{Some nonabelian radical} & \casecell{\(N\le2\nu\)} & \casecell{The centralizer-cover estimate in Section~\ref{sec:binary-nonabelian}.} \\
\casecell{All radicals abelian; rank two occurs} & \casecell{All remaining \(N\)} & \casecell{Lemma~\ref{lem:rank-two-terminal}: \(a(P)=N\).} \\
\casecell{Abelian noncentral radical; all ranks \(\ge4\)} & \casecell{\(N\ge2d+5\)} & \casecell{The uniform fibre average.} \\
\casecell{Same} & \casecell{\(N\le2d+4\)} & \casecell{Lemma~\ref{lem:binary-action-cover}; separate fibres, except for \(d=4\), \(N=9,10\), handled by Lemma~\ref{lem:binary-nine-cover}.} \\
\casecell{All radicals central; one component} & \casecell{All remaining \(N\)} & \casecell{Lemma~\ref{lem:spread}; the equality model.} \\
\casecell{All radicals central; several components} & \casecell{\(\nu\ge10\); or \(\nu=4,6,8\)} & \casecell{The amplified clique forces the cyclic-lift range; otherwise Lemma~\ref{lem:small-binary-cover}.} \\
\hline
\end{tabular}
\end{center}

For the first row, cyclic lifts give \(a(P)\le2^\nu-1<B(N)\)
when \(N\ge9\). If \(N\le8\), then \(\nu\le3\): the cost
is at most three for \(\nu\le2\), and seven for \(\nu=3\),
which requires \(N\ge7\).

In the remaining rows, a nonabelian radical raises every fibre
width; an abelian noncentral radical raises widths in directions
that act nontrivially. With central radicals, we use the
commutator forms themselves.

\subsection{A nonabelian radical}\label{sec:binary-nonabelian}
Suppose some scalar radical \(R\) is nonabelian. Put
\(A=\Z(R)\) and write \([R:A]=2^m\) and
\([A:\Z(P)]=2^s\). Thus \(\nu=d+m+s\), with
\(d=\dim(P/R)\ge2\) and \(m\ge2\).

The radical can be covered by \(2^m-1\) cyclic lifts,
and one such lift gives an abelian subgroup of index
at most \(2^{m-1}\) in \(R\). Thus the fibre bound gives
\[
 a(P)\le 2^m-1+2^{m-1}\frac{(2^d-1)N}{d+1}.
\]
By Lemma~\ref{lem:radical-width}, the width of \(gR\) is at least
\(m+1+\log_2[A:C_A(g)]\). If \(A=\Z(P)\), the widths over a
simplex give \(N\ge(d+1)(m+1)\ge2(d+m)+1=2\nu+1\), since
\(dm\ge d+m\). This contradicts \(N\le2\nu\).

Now \(A\ne\Z(P)\). Conjugation on \(A\) factors through \(P/R\).
The simplex spans \(P/R\), so at least one of its lifts acts
nontrivially on \(A\). Number the lifts \(g_1,\ldots,g_{d+1}\)
so that \(g_{d+1}\) is such a lift. The simplex vertices sum
to zero, so the first \(d\) still span \(P/R\). Thus
\(\bigcap_{i=1}^{d}C_A(g_i)=C_A(P)=\Z(P)\). A coset of this
intersection is determined by its cosets modulo the
\(C_A(g_i)\), giving
\[
 \prod_{i=1}^{d}[A:C_A(g_i)]\ge[A:\Z(P)]=2^s.
\]
The first \(d\) fibres therefore contribute at least \(s\)
vertices beyond their baseline widths, and the last fibre
contributes at least one more. Hence
\[
 N\ge(d+1)(m+1)+s+1=\nu+dm+2.
\]

Together with \(N\le2\nu\), this gives
\(N\ge2dm+4\ge4(d+m-1)\), using \((d-2)(m-2)\ge0\).

Put \(h=d+m-1\) and \(\kappa=\lfloor(N-1)/2\rfloor\).
The clique bound gives \(\kappa\ge2h-1\). The cover estimate
and \((2^d-1)/(d+1)\le2^{d-2}\) give
\[
 a(P)\le 2^{h-2}(N+2)-1\le 2^{h-1}(\kappa+2)-1.
\]
It remains to compare \(\kappa+2\) with
\(2^{\kappa-h+1}\). At \(\kappa=2h-1\), the required
inequality is \(2h+1\le2^h\), which holds since \(h\ge3\).
Thereafter \(2^{\kappa-h+1}\) doubles
with each unit of \(\kappa\), whereas \(\kappa+2\) grows by one.
Thus \(a(P)<B(N)\).

\subsection{An abelian noncentral radical}
Assume now that all scalar radicals are abelian.
If some scalar form has rank two, we have \(a(P)=N\)
by Lemma~\ref{lem:rank-two-terminal}.
For \(N\ge9\), this is strictly below \(T(N)\).
If all scalar ranks are at least four, choose a
noncentral radical of rank \(d\), if one exists, and use
\eqref{eq:uniform-fibre-cost} with \(a(R)=I=1\).
Write \(\kappa=\lfloor(N-1)/2\rfloor\), so the
fibre cost divided by \(B(N)\) is at most
\[
 \frac{2\kappa+2}{d+1}\,2^{d-\kappa}.
\]
For \(\kappa\ge d+2\) the ratio decreases with \(\kappa\), starting
at \((d+3)/(2(d+1))\le7/10\).
When \(\kappa\le d+1\), we have \(N\le2d+4\). The averaging
bound alone is insufficient here, since it does not exclude
individual fibres of large width. We first show that every
fibre has width at most two. Covering the fibres separately
then suffices except at \(d=4\)
with \(N=9,10\), where we build a cover of the central
quotient by common isotropic subspaces.

\begin{lemma}[Width reduction in small scalar rank]\label{lem:binary-action-cover}
Let \(R\) be an abelian noncentral scalar radical in a finite
\(2\)-group \(P\), with \(d=\dim(P/R)\ge4\). Assume all nonzero
scalar forms have rank at least four. If \(N=\omega(P)\le2d+4\),
then \(d\in\{4,6\}\), \(N\ge9\), and every fibre of \(P/R\)
has width at most two.
\end{lemma}

\begin{proof}
\emph{Step 1: cliques from two fibres.}
Write \(R\) additively and set \(D_g=\alpha_g-1\), where
\(\alpha_g\) is right conjugation on \(R\).
Since \(R\) is abelian, this action factors through the elementary
abelian group \(V=P/R\). In particular the actions commute,
and the notation \(D_v\), for \(v\in V\), is unambiguous.
The width of the fibre at \(v\) is \(|D_v(R)|\), by
Lemma~\ref{lem:abelian-coset}. The commuting
classes in \(gR\) are the cosets of \(K_g=\ker D_g\).
For distinct fibres \(yR,zR\), two elements \(yr,zs\)
commute exactly when
\begin{equation}\label{eq:binary-two-fibre}
 D_z(r)-D_y(s)=-[y,z].
\end{equation}
Indeed, \(yrzs=yz\,r^zs\) and
\(zsyr=yz[y,z]^{-1}s^yr\); all the factors after \(yz\)
belong to the abelian subgroup \(R\).
Suppose both fibres have width two and their kernels
differ. Start with a noncommuting pair. Vary \(r\) by
an element of \(K_z\setminus K_y\) and \(s\) by an element
of \(K_y\setminus K_z\). Neither change affects the
equation, and each gives a second class in its own
fibre. This yields a four-point clique.

More generally, fibres with different actions give at
least three clique points. A fibre of width at least
four suffices on its own. If one fibre has width one and
the other width two, choose both classes in the latter.
Equation~\eqref{eq:binary-two-fibre} lets us choose a point in the inactive
fibre that commutes with neither. If both widths are two
and the kernels agree, the two operators have distinct
images of order two. The equation has at most one solution
in the pair of class labels, so we can use both classes
in one fibre and one class in the other.

\emph{Step 2: restricting the scalar dimension.}
Choose a direction \(x\in V\) acting nontrivially.
The \(d-1\) points of a symplectic simplex in
\(x^\perp/\langle x\rangle\) each lift to two directions
\(y,y+x\). Their actions differ, since
\(\alpha_{y+x}=\alpha_y\alpha_x\) and \(\alpha_x\ne1\).
Each pair supplies at least three points by Step 1,
and the scalar form separates different pairs. Thus
\(3(d-1)\le N\le2d+4\), which leaves \(d=4,6\) and \(N\ge9\).

\emph{Step 3: linearizing the action.}
We next show that every fibre has width at most two. The
first step is to make the action a vector space of maps.
If the directions of width at most two lay in a hyperplane,
the simplex in Lemma~\ref{lem:clique-meets-hyperplane}
would give at least \(4d+1\) vertices. This contradicts \(N\le2d+4\), so those
directions span \(V=P/R\). Their action operators have
images of order at most two. Every other action preserves
each such image because the actions commute. An image of
order two has no nonidentity automorphism, so it is fixed
pointwise. Thus \(D_uD_v=0\) for the spanning directions, and
\(D_{u+v}=D_u+D_v+D_uD_v\) makes all operators sums of the spanning
ones. Their images generate the central elementary abelian
group \([R,P]\), and \(R^2\le\Z(P)\).
The common kernel in \(R\) of all the operators is
\(C_R(P)=\Z(P)\), since the centre lies in every scalar radical.
Thus \(R/\Z(P)\) and \([R,P]\) are binary vector spaces.
The assignment \(v\mapsto D_v\) is now linear, and its image
\(\mathcal D\le\operatorname{Hom}_{\F_2}(R/\Z(P),[R,P])\)
is spanned by maps of rank one. A fibre has width
\(2^{\operatorname{rank}D_v}\).

\emph{Step 4: excluding wider fibres except at one endpoint.}
If some direction \(x\) has \(\operatorname{rank}D_x\ge2\),
repeat the clique construction in \(x^\perp/\langle x\rangle\).
Each pair of lifts \(y,y+x\) now supplies four vertices.
Indeed, if both widths are at most two, the maps
\(D_y,D_{y+x}\) must have rank one and distinct kernels,
since their sum has rank at least two. The two-fibre
argument applies. Hence \(N\ge4(d-1)\). In the range
\(N\le2d+4\), this leaves only \(d=4,\ N=12\).

\emph{Step 5: excluding the endpoint \(d=4\), \(N=12\).}
At this endpoint, we need two rank-one maps \(A_1,A_2\)
with different images and different kernels. If no pair
differed in both respects, the rank-one maps would share
an image line or a kernel. Indeed, two with distinct images
would have the same kernel, and every other rank-one map
would have to share it. Their sums would then all have
rank at most one, contrary to the existence of a map
of rank at least two.

In the four-space \(V\), take the plane spanned by two
directions inducing \(A_1,A_2\). Each translate of the plane
contains a direction of width at least four. To see this,
choose inputs \(u,v\) with \(A_1(u)=a\ne0\), \(A_2(v)=b\ne0\)
and \(A_2(u)=A_1(v)=0\). The outputs \(a,b\) are independent.
In each translate, the two possible values at \(u\)
differ by \(a\), and the two at \(v\) differ by \(b\).
The choices are independent. If every pair were dependent,
any nonzero value at \(v\) would span a line containing
both values at \(u\), and hence containing \(a\).
Both values at \(v\) would then lie on \(\langle a\rangle\),
contrary to their difference \(b\).

The four translates give at least four wide directions.
A totally isotropic subspace in the symplectic four-space
has at most three nonzero vectors. Two wide directions
\(x,y\) therefore pair nontrivially. Put
\(U=\langle x,y\rangle^\perp\). We complete \(x,y\) to
a simplex with the three directions \(x+y+u\) for
\(0\ne u\in U\). They pair nontrivially with \(x,y\)
and with each other.

At least two of these last three directions act nontrivially.
If two acted trivially, the kernel of the action map
\(V\to\mathcal D\) would have dimension at least two.
Rank-nullity then gives \(\dim\mathcal D\le2\); equality
holds since \(\mathcal D\) contains \(A_1,A_2\).
Its only rank-two map would be \(A_1+A_2\), so \(D_x=D_y\).
The actions at \(x+y+u\) would then be \(D_u\).
Two trivial actions would force all of \(U\) to act
trivially. But then \(V\) would act through just the
one map \(D_x\), contradicting \(\dim\mathcal D=2\).
The simplex therefore gives at least \(4+4+2+2+1=13\)
clique vertices.

This excludes rank at least two throughout the range
under consideration. Every fibre has width at most two.
\end{proof}

Covering each fibre separately now costs at most \(2^{d+1}-1\).
For \(N\ge2d+3\), this is below \(B(N)\ge2^{d+1}\).
Only \(d=4\) and \(N=9,10\) remain. Here we replace the
separate covers by subspaces isotropic for every commutator form.

\begin{lemma}[Nine common isotropic subspaces]\label{lem:binary-nine-cover}
Let \(R\) be an abelian noncentral scalar radical in a finite
\(2\)-group \(P\), with \([P:R]=16\). Suppose every nonzero
scalar form has rank at least four, every fibre has width at
most two, and \(N=\omega(P)\le10\). Then \(P/\Z(P)\cong\F_2^5\)
has a cover by nine common isotropic subspaces. Consequently
\(a(P)=N=9<T(N)\).
\end{lemma}

\begin{proof}
\emph{Step 1: the five-dimensional quotient.}
We first prove that \([R:\Z(P)]=2\) and that \(P\) has
class two with central squares. This reduces the problem to
an alternating map on a five-dimensional binary space.
We then construct the cover from a Kronecker block when
there are two scalar forms, and from the incidence of common
isotropic planes when there are three.

Put \(V=P/R\) and retain the action notation \(D_g=\alpha_g-1\)
from the preceding proof. These operators form a space of
maps from \(R/\Z(P)\) to \([R,P]\),
each of rank at most one.
First we show \([R:\Z(P)]=2\); a larger quotient would give
an eleven-point clique.
If \(\dim(R/\Z(P))\ge2\), the kernels of the rank-one maps
cannot all coincide, since they intersect in zero. The maps must
therefore share an image line. Here we use the following
rank-one observation: if \(a\otimes f\) and \(b\otimes g\)
have different image lines and different kernels, their sum
has rank two. Thus any two rank-one maps in \(\mathcal D\)
share an image or a kernel. Once two have different kernels,
they share an image, and comparison with both shows that
every other nonzero map has that same image.
Their action functionals
span the dual of \(R/\Z(P)\). Project these functionals onto
two coordinates, giving a surjection \(\lambda:V\to\F_2^2\)
with kernel \(J\).

Choose \(x\notin J\cup J^\perp\). Then \(\lambda(x)\ne0\)
and \(x^\perp\cap J\) has dimension one, so
\(\lambda\) maps \(x^\perp\) onto \(\F_2^2\).
Write \(00,10,01,11\) for the vectors of \(\F_2^2\), and
choose coordinates with \(\lambda(x)=10\). The second
coordinate is then a nonzero functional on the
symplectic plane \(x^\perp/\langle x\rangle\).
It is one at two of the three nonzero points and zero
at the third. In each pair of lifts, adding \(x\) toggles
the first coordinate. The two pairs over the points
with second coordinate one thus carry the labels
\(01,11\), and the remaining pair carries \(00,10\).

In the first two pairs, the action functionals are
nonzero and distinct because their projected labels are.
Over \(\F_2\) they therefore have distinct kernels. Each pair gives four clique points by the two-fibre argument
in Lemma~\ref{lem:binary-action-cover}. In the last pair the
actions differ, so the same lemma gives at least three points. The three nonzero
points of a symplectic plane pair nontrivially, so
the \(4+4+3\) points form a clique, contrary to \(N\le10\).

Thus \([R:\Z(P)]=2\). We show that the central quotient
is an elementary abelian space of dimension five, so
we can look for a cover by common isotropic subspaces.
The action of \(P\) on the group \(R/\Z(P)\) of order two
is trivial.
If \(P'\not\le\Z(P)\), projection of the commutator onto
\(R/\Z(P)\cong\F_2\) gives a nonzero scalar form on \(P/R\).
It is well defined because \([R,P]\le\Z(P)\). Moreover
\(\gamma_3(P)\le[R,P]\le\Z(P)\), so the projection kills
\(\gamma_3(P)\): this is one of the scalar forms to which
the rank hypothesis applies. On the four-space \(P/R\),
it is therefore nondegenerate.
Choose \(t\) acting nontrivially on \(R\). Its perpendicular
for this form contains two vectors with nonzero pairing,
by rank at least four. For lifts \(x,y\) of these vectors,
conjugation by \(t\) changes \(x,y\) by central factors,
so it fixes \([x,y]\). Yet \([x,y]\) represents the nonzero
class in \(R/\Z(P)\), whose representatives are moved by
\(t\). This proves that \(P\) has class two.

Its squares lie in \(R\). Since \([x^2,y]=[x,y^2]\),
the commutator \([x^2,y]\) is trivial if either square
is central. If both are noncentral, the squares agree
modulo \(\Z(P)\), so \(x^2\) again commutes with \(y\).
Thus every square is central, and \([x,y]^2=[x^2,y]=1\)
shows that \(P'\) is elementary abelian.

\emph{Step 2: the alternating map on the quotient.}
Write \(E=P/\Z(P)\cong\F_2^5\), \(W=P'\) and let \(\beta\)
be the commutator map. Write \(C_\beta(x)=\ker\beta(x,\cdot)\)
for the common centralizer of \(x\) in \(E\).
The original scalar form vanishes
on the nonzero space \([R,P]\), so \(\dim W\ge2\).
All nonzero scalar forms have rank four.

\emph{Step 3: two scalar forms.}
If \(\dim W=2\), Scharlau's decomposition \cite{scharlau}
gives a single singular Kronecker block of dimension five.
A singular block of dimension one would give a nonzero
common radical. If a singular block had dimension three,
the remaining dimension two could contain neither another
singular block (there is no common radical) nor a regular
block nondegenerate for all three nonzero combinations:
the space of alternating forms on a two-space has dimension
one. Some nonzero combination would consequently have rank
at most two, contrary to the rank-four hypothesis.
The only possible decomposition is therefore the single
five-dimensional singular block. A singular Kronecker block
of a prescribed odd dimension is unique up to simultaneous
equivalence of the two alternating forms \cite{scharlau}.
It is therefore enough to construct the cover in any model
of this block. To exhibit such a model,
take \(K=\F_8\), a two-dimensional subspace \(U\le K\), and
a surjection \(\pi:K\to\F_2^2\). On \(K\oplus U\) the pair
\[
 \beta((x,s),(y,t))=\pi(xt-ys)
\]
has zero common radical. Multiplication by \(s\ne0\) is
onto \(K\), and \(xU\) cannot lie in \(\ker\pi\) for \(x\ne0\).
Each nonzero scalarization has rank four because the pairing
between \(K\) and \(U\) has rank two. It therefore models
the same Kronecker block. The space \(K\oplus0\) and the
eight graphs \(\{(as,s):s\in U\}\), \(a\in K\), partition
the nonzero vectors into common isotropic subspaces.

\emph{Step 4: at most three scalar forms.}
Suppose \(\dim W\ge3\). We look for a three-space
among the radical vectors of the scalar forms. There is
one nonzero radical vector per nonzero scalar form, since
every form has rank four.

If there were four independent forms, the fifteen nonzero
combinations would have at most fifteen radical vectors.
At a nonzero vector \(x\) outside that set, the four
functionals obtained by evaluation at \(x\) are independent:
a linear relation would put \(x\) in the radical of a nonzero
combination. Their common kernel has dimension \(5-4=1\)
and contains \(x\), so it is exactly \(\langle x\rangle\).
Two distinct vectors outside the set therefore cannot
commute even for these four forms. There are at least
\(31-15=16\) such vectors, contrary to \(N\le10\).
Thus \(\dim W=3\).

\emph{Step 5: counting the common isotropic planes.}
Let \(S\) be the set of radical vectors of the seven
nonzero forms. Outside \(S\cup\{0\}\), evaluation has
rank three, so the common centralizer is a plane.
Each such vector lies on just that one common isotropic
plane: an isotropic plane containing \(x\) must lie in
\(C_\beta(x)\), which already has dimension two.
Take one vector outside \(S\) from each plane that meets
the complement of \(S\cup\{0\}\). If two representatives
commuted, they would span the unique isotropic plane through
either of them. They would therefore come from the same
plane. The selected representatives thus form a clique.

To see how many planes this leaves inside \(S\), count
the planes through a nonzero vector \(x\). If \(n_x\)
nonzero forms have radical \(\langle x\rangle\), the
annihilator of \(\beta(x,E)\) has size \(n_x+1\).
Rank-nullity gives \(|\ker\beta(x,\cdot)|=4(n_x+1)\).
The planes through \(x\) correspond to the nonzero
vectors of \(\ker\beta(x,\cdot)/\langle x\rangle\), so
there are
\[
 \frac{4(n_x+1)}2-1=1+2n_x.
\]
Each of the seven nonzero scalar forms has exactly one
nonzero radical vector, so \(\sum_{x\ne0}n_x=7\), even
if several forms have the same radical. Summing over
\(x\ne0\) gives \(31+2\cdot7=45\) incidences.
Each plane has three nonzero vectors, hence there are
exactly fifteen common isotropic planes. The clique bound
\(N\le10\) leaves at least five planes with all their
nonzero vectors in \(S\).

\emph{Step 6: forcing an isotropic three-space.}
These five planes force the three-space. Two must meet,
since \(S\) has at most seven points. If they meet at
\(x\ne0\), their span \(Y\) has dimension three and
already contains five vectors of \(S\). A plane not
contained in \(Y\) uses at least two vectors outside
\(Y\), and there are at most two available in \(S\).
Those two vectors determine the plane, so at most one
of the five planes lies outside \(Y\). Among the at least
four planes in \(Y\), one avoids \(x\), since only three
can pass through \(x\). It is isotropic, and \(x\)
commutes with all of \(Y\) by the construction of \(Y\).
Together this plane and \(x\) span a common isotropic
three-space.

\emph{Step 7: partitioning the remaining vectors.}
For \(0\ne y\in Y\), we have \(Y\le C_\beta(y)\), so \(y\in S\).
These seven vectors exhaust \(S\). Since each accounts for
one of the seven scalar forms, \(C_\beta(y)=Y\) for \(0\ne y\in Y\).
Every vector outside \(Y\) lies on a unique common isotropic
plane, and this plane meets \(Y\) only at zero. These planes
partition the \(24\) vectors outside \(Y\) into eight triples,
giving eight planes in addition to \(Y\).

In either case the nine subspaces lift to nine abelian
subgroups. Lemma~\ref{lem:binary-action-cover} applies to the
present hypotheses and gives \(N\ge9\). Together with
\(N\le a(P)\le9\), this yields \(a(P)=N=9<T(N)\).
\end{proof}

\subsection{Central radicals}
If every scalar radical is central, then \(P\) has class
two, \(P/\Z(P)\) is elementary abelian of even dimension
\(\nu\), and every nonzero scalarization of the commutator
is nondegenerate. If the commutator has one component,
the spread gives \(a(P)=2^{\nu/2}+1\), with \(N=\nu+1\).
For \(\nu=2,4,6\) these cover sizes are \(3,5,9\)
and fit the stated profile. For \(\nu\ge8\), this
group has the binary symplectic model and \(a(P)=T(N)\).

With several forms, we again use one to separate fibres
and a second to enlarge the clique inside them. Here we
can triple each contribution by choosing a plane on which
the forms behave differently. Take independent forms
\(b,c\). There is a plane
\(J\) isotropic for \(b\) and nondegenerate for \(c\).
To find it, write \(c(x,y)=b(Sx,y)\), choose \(x\)
with \(x,Sx\) independent, and choose
\(y\in x^{\perp_b}\) with \(c(x,y)\ne0\).
Then \(J=\langle x,y\rangle\) has the desired properties.
In \(J^{\perp_b}/J\), choose a binary symplectic clique of
size \(\nu-3\), using the singleton \(\{0\}\) if \(\nu=4\).
Lift each vertex to a representative
orthogonal to \(J\) for \(c\), which is possible because
\(c|_J\) is nondegenerate. Adding the three nonzero
vectors of \(J\) to this representative gives three
pairwise \(c\)-nonorthogonal lifts. Lifts of different
vertices are already separated by \(b\). Thus
\begin{equation}\label{eq:binary-central-clique}
 N\ge3(\nu-3).
\end{equation}
When \(\nu\ge10\), this implies \(N\ge2\nu+1\), so
cyclic lifts suffice. In the smaller dimensions, we
use subspaces isotropic for the full commutator. In
dimensions four and six, spreads or common centralizers
supply them. In dimension eight, we first cover disjoint common
isotropic planes, then pair commuting points among the remaining vectors.
The planes reduce the number of subgroups enough to obtain
a strict bound.

\begin{lemma}[Regular pencils in dimensions four and six]\label{lem:regular-pencil}
Let \(b,c\) be independent alternating forms on a
\(2j\)-dimensional space \(V\) over \(\F_p\), with \(j=2\)
or \(3\). If every nonzero form in their span is
nondegenerate, then \(V\) can be identified with
\(\F_{p^j}^2\) so that its \(\F_{p^j}\)-lines are
isotropic for both forms.
\end{lemma}

\begin{proof}
In Scharlau's regular decomposition \cite{scharlau},
a block with elementary divisor \(f^k\) has dimension
\(2k\deg f\). Here there are no singular blocks or linear
elementary divisors, since every nonzero combination is
nondegenerate. In particular each irreducible \(f\) has
degree at least two. The sum of the quantities \(k\deg f\)
over the regular blocks is \(j\). Since \(j=2\) or \(3\),
two blocks would already contribute at least four, and
\(k\ge2\) in one block would also contribute at least four.
Thus there is exactly one block, with \(k=1\) and
\(\deg f=j\). On this block the two forms
are scalarizations of the determinant on \(\F_{p^j}^2\).
The \(\F_{p^j}\)-lines are therefore isotropic for both.
\end{proof}

\begin{lemma}[Small binary central quotients]\label{lem:small-binary-cover}
Let \(P\) have class two and elementary abelian central
quotient of dimension \(\nu\in\{4,6,8\}\). Suppose that
there are at least two independent scalar commutator
forms, and that every nonzero form is nondegenerate.
Then \(a(P)\le T(N)\), where \(N=\omega(P)\).
For \(N\ge9\) the inequality is strict.
\end{lemma}

\begin{proof}
\emph{Step 1: the common-centralizer dimension.}
Write \(V=P/\Z(P)\), let \(\beta\) be the commutator
map on \(V\), and put \(N=\omega(P)\) and
\(e=\dim_{\F_2}P'\). The Pfaffian bound \cite[Lemma~3]{gow} gives
\(e\le\nu/2\).

For every \(z\ne0\), the evaluation map
\(\beta(z,\cdot):V\to P'\) is onto: otherwise a nonzero
functional annihilating its image would give a scalar form
with \(z\) in its radical. Hence every common centralizer
has dimension \(\nu-e\).

\emph{Step 2: dimension four.}
If \(\nu=4\), then \(e=2\), and the regular pencil gives
five blocks by Lemma~\ref{lem:regular-pencil}. The scalar
clique has size five, so \(a(P)=N=5\).

\emph{Step 3: dimension six, a common spread or a third form.}
If \(\nu=6\), then \(e=2\) or \(3\) and
\eqref{eq:binary-central-clique} gives \(N\ge9\).
With two forms, the regular pencil gives nine blocks.
With three forms and \(N\ge11\), first apply the regular
pencil lemma to two forms, giving nine common isotropic
three-spaces. On each three-space the restriction of the
third alternating form has rank zero or two. In rank zero
keep the whole space. In rank two its radical is a line,
and the inverse images of the three lines of the
two-dimensional quotient are isotropic planes covering the
three-space. Thus each of the nine blocks costs at most three.
This costs at most \(27<B(N)\).

\emph{Step 4: dimension six, the range \(N\le10\).}
With three forms, put \(Y_z=\ker\beta(z,\cdot)\). Each \(Y_z\)
has dimension three by the evaluation calculation above.
We show that all these spaces are common isotropic, unless
there is an eleven-point clique. If all \(Y_z\) are common isotropic, and \(0\ne t\in Y_z\),
then \(Y_z\subseteq Y_t\). Both have dimension three, so
\(Y_z=Y_t\). Thus two such spaces meeting off zero coincide.
They partition the \(63\) nonzero vectors into blocks of
seven, giving nine blocks.

Otherwise choose \(z\) with \(Y=Y_z\) nonisotropic and a
scalar form \(c\) nonzero on \(Y\). In the symplectic space
\(z^{\perp_c}/\langle z\rangle\), the plane
\(J=Y/\langle z\rangle\) is nondegenerate. Translate the
three nonzero vectors of \(J\) by a nonzero \(u\in J^\perp\)
and add the two nonzero vectors of \(J^\perp\) other than
\(u\). The five vectors lie outside \(J\). Within each set the
pairings are unchanged, and across the sets they are the
nonzero pairings with \(u\). Lift each vector to its two translates
by \(z\). The form \(c\) separates lifts of different
vectors. If \(v\) is one of the chosen lifts, its image lies
outside \(J\), so \(v\notin Y\) and \(\beta(v,z)\ne0\).
Hence \(v\) and \(v+z\) fail to commute with one another,
and both fail to commute with \(z\). The five doubled
directions together with \(z\) therefore give eleven vertices,
a contradiction.

\emph{Step 5: dimension eight.}
Here \(N\ge15\) by \eqref{eq:binary-central-clique}.
If \(N\ge17\), cyclic lifts suffice since \(z(P)=256\le B(N)\).
For \(N=15,16\), pairing commuting points could leave up
to \(N\) singletons. We first choose \(N\) disjoint common
isotropic planes. Each covers three points in one subgroup,
offsetting the cost of one singleton.

Every nonzero vector lies on the same positive number \(D\)
of common isotropic planes. Indeed, \(\beta(x,\cdot)\)
is onto \(P'\) for \(x\ne0\), so those planes are the lines
in \(\ker\beta(x,\cdot)/\langle x\rangle\), of dimension
\(7-e\ge3\). There are \(255D/3=85D\) planes in total,
and each chosen plane meets at most \(3D\) planes off zero,
including itself. After \(k<N\le16\) choices, at most
\(3kD\le45D<85D\) planes are forbidden. A greedy choice
therefore gives \(N\) planes with pairwise disjoint sets of
nonzero vectors.
Cover their \(3N\) points in \(N\) subgroups, then pair
commuting points among the rest until no pair remains.
The unpaired points form a clique and hence number at most
\(N\). The total cost is at most
\[
 N+\left\lfloor\frac{255-3N+N}{2}\right\rfloor
 =127<B(N).\qedhere
\]
\end{proof}

\subsection{Products of Sylow subgroups}

An abelian group needs one subgroup. If only one Sylow
subgroup is nonabelian, the abelian factors preserve \(a\), \(\omega\) and the isoclinism class.
Otherwise split \(G\) into two nonabelian direct factors.
Their clique numbers \(i,j\ge3\) satisfy \(ij\le N\),
and multiplying the two covers gives an upper bound for
the product. We use \(T(k)\le2^{(k+1)/2}\) for \(k\ge3\),
direct from the formula for \(T\) (for \(3\le k\le6\),
use \(k^2\le2^{k+1}\)). Unless \(i=j=3\), we have
\((i-1)(j-1)\ge6\), so
\[
 T(i)T(j)\le2^{(i+j+2)/2}
 \le2^{(ij-3)/2}\le B(ij)/\sqrt2.
\]
For \(i=j=3\) the cost is nine. Induction on the number of
nonabelian factors gives \(a(G)\le B(N)<T(N)\), since
\(B\) is nondecreasing.

\section{Abelian covers of nonnilpotent groups}\label{sec:nonnil}

\begin{proposition}[Nonnilpotent groups]\label{prop:nonnil-cover}
If \(G\) is nonnilpotent and \(N=\omega(G)\), then
\(a(G)\le T(N)\), with strict inequality for \(N\ge9\).
\end{proposition}

\paragraph{Case map.}
We argue by strong induction on \(N\). The nilpotent case is
already available, as is the profile for smaller clique numbers.
\begin{center}
\renewcommand{\arraystretch}{1.12}
\begin{tabular}{@{}p{.28\textwidth}p{.25\textwidth}p{\dimexpr0.47\textwidth-4\tabcolsep\relax}@{}}
\hline
\casecell{Group or quotient} & \casecell{Clique range} & \casecell{Covering argument} \\
\hline
\casecell{Any nonnilpotent group} & \casecell{\(N\ge45\)} & \casecell{Centralizer descent and induction on the smaller clique number.} \\
\casecell{Nonsolvable group} & \casecell{\(N\le44\)} & \casecell{A central extension of \(A_5\) of index at most two.} \\
\casecell{Affine, cyclic complement; \(H\) abelian} & \casecell{\(s=1\); or \(s\ge2\)} & \casecell{Cover the nilpotent translation preimage; otherwise the quadratic bound.} \\
\casecell{Affine, cyclic complement; \(H\) nonabelian} & \casecell{\(N\ge(m+1)r\), \(r\ge3\)} & \casecell{Lemma~\ref{lem:affine-index} and the weighted fibre estimate.} \\
\casecell{Affine, nonabelian complement} & \casecell{\(21\le N\le44\)} & \casecell{Affine degree nine and a cover from the widths of translation conjugates.} \\
\hline
\end{tabular}
\end{center}
Here \(m\) is the affine degree, \(s\) a complement-fibre
width, and \(r\) a selected fibre width, as defined below.

\subsection{Centralizer descent for large clique numbers}
\label{sec:centralizer-descent}

The selection argument is the one in
\cite[Proposition~4.3]{index}; we give the proof because it
uses no central-index bound. We use the classical
Baer--Suzuki theorem \cite[Theorem~3.8.2]{gorenstein}:
if \(X\) is a \(p\)-subgroup and \(\langle X,X^g\rangle\)
is a \(p\)-group for every \(g\in G\), then \(X\le O_p(G)\).
Its nilpotent form says that if every two conjugates of a
cyclic subgroup generate a nilpotent group, its normal
closure is nilpotent. Indeed, apply the theorem to the
\(p\)-part of a generator for each prime \(p\); each such
part lies in \(O_p(G)\), so the generator lies in \(F(G)\).
See also \cite[pp.~281--282]{baersuzuki} for this formulation;
the stronger extensions proved there using the classification
of finite simple groups are not needed.

We first record two elementary facts used in the descent.
A nonnilpotent group has clique number at least four.
Indeed, if \(\omega(U)\le3\), \eqref{eq:tomkinson-coset}
with the whole group as its coset gives \(z(U)\le6\).
A central quotient cannot be nontrivial cyclic, so its only
nontrivial possibilities are \(C_2^2\) and \(S_3\).
The latter has a four-point clique (the three transpositions
and a \(3\)-cycle), which lifts to \(U\). The former makes
\(U\) nilpotent of class at most two.

Also, if \(K\) and \(L\) centralize one another, then
\begin{equation}\label{eq:commuting-factor-cliques}
 \omega(KL)\ge\omega(K)\omega(L/\Z(L)).
\end{equation}
To see this, multiply every vertex of a clique in \(K\) by
every chosen lift of a clique in \(L/\Z(L)\).
For distinct quotient vertices, a commutator from \(L\)
is noncentral in \(L\), whereas any cancelling commutator
from \(K\) would lie in \(K\cap L\le\Z(L)\).
For equal quotient vertices, noncommutation is inherited
from \(K\). The products are distinct for the same reason.

\begin{lemma}[Centralizer descent]\label{lem:centralizer-descent}
Let \(G\) be a finite nonnilpotent group and put \(N=\omega(G)\).
There is a nonnilpotent subgroup \(L\le G\) whose centralizer
\(K=C_G(L)\) satisfies
\[
 [G:K]\le N(N-1),\qquad 4\omega(K)\le N.
\]
\end{lemma}

\begin{proof}
Put \(F=F(G)\), \(S=G\setminus F\), and \(f=[G:F]\).
Take a maximum clique \(D\) in \(S\), of size \(t\le N\).
Its centralizers cover \(S\). For \(C=C_G(x)\) with \(x\in D\),
we have \([C:C\cap F]\le f\), so
\(|C\cap S|\le(1-1/f)|C|\). Since \(|S|=(1-1/f)|G|\), summing gives
\[
 \sum_{x\in D}\frac{1}{[G:C_G(x)]}\ge1.
\]
Some \(x\in D\) therefore has \(r=[G:C_G(x)]\le t\le N\).
By Baer--Suzuki, \(x\) has a conjugate \(y\) such that
\(L=\langle x,y\rangle\) is nonnilpotent. The \(C_G(x)\)-orbit
of \(y\) avoids \(x\), so it has size at most \(r-1\). Its
stabilizer is \(K=C_G(x,y)\), giving \([G:K]\le r(r-1)\).
The quotient \(L/\Z(L)\) is nonnilpotent, since a central
extension of a nilpotent group is nilpotent. It has a clique
of size four. Now \eqref{eq:commuting-factor-cliques} gives
\(4\omega(K)\le N\).
\end{proof}

\paragraph{The numerical cost of descent.}
Let \(N\ge45\) and put \(j=\lfloor N/4\rfloor\).
For the centralizer in Lemma~\ref{lem:centralizer-descent},
induction gives \(a(K)\le T(j)\), and \eqref{eq:lift} gives
\[
 a(G)\le T(j)\{1+(N(N-1)-1)(N-1)\}<N^3T(j).
\]
For \(45\le N\le48\) this bound is at most
\(48^3\cdot33<2^{22}=B(45)\).
Increasing \(N\) by four quadruples \(B(N)\), while \(T(j)\)
grows by at most a factor two and \(N^3\) grows by a factor
at most \((49/45)^3<2\). Thus \(a(G)<B(N)<T(N)\) for
\(N\ge45\).

\subsection{Nonsolvable groups with \texorpdfstring{\(N\le44\)}{N at most 44}}

Let \(G\) be nonsolvable with \(N\le44\), and let \(S\) be
its solvable radical. Write \(G''=(G')'\). In alternatives
(1) and (2) of \cite[Theorem~1.1]{nonsol}, \(G=\Z(G)K\)
with \(K\cong\mathrm{PSL}_2(p)\) or \(\mathrm{SL}_2(p)\),
where \(p\in\{5,7\}\). Here \(G''=K\), since \(K\) is perfect,
and \(S=\Z(G)\), since \(G/\Z(G)\) is nonabelian simple.
Taking \(a=1\) in these two cases, all four alternatives give
\(G=G''\langle a\rangle S\) with \(a^2\in\Z(G)\) and \(G''\)
isomorphic to \(A_5,\mathrm{SL}_2(5),\mathrm{PSL}_2(7)\) or
\(\mathrm{SL}_2(7)\). The last two are excluded by
\(\omega(G''/\Z(G''))\le N\) and
\(\omega(\mathrm{PSL}_2(7))=57\) \cite[Lemma~4.4]{noncommuting-graph}.
Thus \(G''/\Z(G'')\cong A_5\).

The nonidentity elements of \(A_5\) have as centralizers its five
Klein four subgroups, ten subgroups of order three and six
subgroups of order five. Two commuting nonidentity elements
have the same centralizer, since these centralizers are abelian.
Choosing one element from each gives \(\omega(A_5)\ge21\).

The intersection \(G''\cap S\) is solvable and normal in \(G''\),
so its image in \(A_5\) is trivial. Conversely \(\Z(G'')\) is
abelian and normal in \(G\), hence \(G''\cap S=\Z(G'')\).
In particular \([G'',S]\le\Z(G'')\). For each \(s\in S\),
the map \(x\mapsto x^{-1}x^s\) from \(G''\) to its centre is
a homomorphism. It is trivial by perfectness, so \([G'',S]=1\).
By \eqref{eq:commuting-factor-cliques},
\(21\omega(S)\le\omega(G''S)\le44\), forcing \(S\) to be abelian.

Put \(M=G''S\). Then \(S\le\Z(M)\) and \(M/S\cong A_5\).
The \(15+10+6=31\) nontrivial cyclic subgroups of \(A_5\)
have abelian preimages in \(M\) and cover it.
Since \(a^2\in\Z(G)\le S\), we have \([G:M]\le2\).
By \eqref{eq:lift},
\[
 a(G)\le31N<B(N)\qquad (N\ge21).
\]
Indeed \(31\cdot22<2^{10}=B(21)=B(22)\), and in steps of
two the left side grows by less than a factor two while
\(B\) doubles.

\subsection{The affine quotient}

Let \(G\) be solvable and nonnilpotent with \(N\le44\).
Choose a nonnilpotent quotient \(Q=G/H\) of smallest order.
Every proper quotient of \(Q\) is nilpotent. Take a nonnormal
maximal subgroup \(L\) of \(Q\).
It has trivial core, since a nontrivial core would
give a smaller quotient with a nonnormal maximal
subgroup, hence a smaller nonnilpotent quotient.

The primitive solvable structure theorem \cite[Theorem~A.15.6]{doerk-hawkes}
gives \(Q=V\rtimes L\), with \(|V|=m=p^d\) and a
faithful irreducible action of \(L\) on \(V\).
We view \((v,a)\in Q\) as the affine map
\(x\mapsto v+a(x)\). Thus \(V\) consists of translations
and \(L\) of maps fixing the origin. The inverse images in \(G\) of affine maps are the fibres
whose widths we will add.

Since \(L\) acts faithfully and \(Q\) is nonnilpotent,
we have \(m\ge3\). The complement \(L\cong Q/V\) is nilpotent by
minimality, and its largest normal \(p\)-subgroup \(O_p(L)\) is trivial.
Thus \(p\nmid|L|\). A clique of \(Q\) lifts to a clique
of \(G\). Any nonidentity central element \(a\in L\)
acts without nonzero fixed points, since its fixed
space is \(L\)-invariant and the action is faithful.
The \(m\) elements with linear part \(a\) form a
clique, to which any nonzero translation can be
adjoined. Hence \(m+1\le N\).

\paragraph{When the complement is abelian.}
The algebra generated by \(L\) is a field by
irreducibility, so \(L\) is cyclic of order \(e>1\)
dividing \(m-1\). If the linear part \(a\) is nonidentity,
then \(1-a\) is invertible, so the affine map has a
unique fixed point. It therefore belongs to exactly one
of the \(m\) point stabilizers, which are the conjugates
of \(L\). An affine map commuting with it must fix that
point. Thus nonidentity maps from different stabilizers
do not commute, and a nonzero translation commutes with
none of them. This lets us add clique sizes across
the corresponding fibres.

Let \(K\) be the preimage of the translations.
Write \(u\) for the largest width of a nonzero translation
fibre and \(s\) for the largest width of a nonidentity
complement fibre. By conjugacy a fibre of width \(s\)
occurs over each point stabilizer. Choose one there and a
translation fibre of width \(u\). The cliques in these
fibres join by Lemma~\ref{lem:fibre-calculus}(ii). Writing
\(v=ms\) for the contribution from the complements, we
get \(u+v\le N\).

By Lemma~\ref{lem:fibre-calculus}(i), select one translation
fibre per line and one generator per nontrivial cyclic
subgroup of each complement. There are \(c=(m-1)/(p-1)\)
translation fibres and \(\ell\) fibres per complement.
In particular \(\ell\le\lfloor e/2\rfloor\), the number
of inversion orbits on nonidentity elements.
For an abelian subgroup of index \(I\) in \(H\),
the quotient estimate gives
\begin{equation}\label{eq:affine-kernel-cost}
 a(G)\le a(H)+I(cu+\ell v).
\end{equation}
If \(m\) is odd, then \(p\ge3\)
and \(c,\ell\le(m-1)/2\). If \(m\) is even, then
\(c\le m-1\) and \(\ell\le(m-2)/2\).

\subsection{Abelian kernel and cyclic complement}

Write \(\kappa=\lfloor(N-1)/2\rfloor\), so \(N\le2\kappa+2\).
For \(\kappa\ge4\) we have \(\kappa^2\le2^\kappa=B(N)\), since the
square starts at \(4^2=2^4\) and grows by less than
a factor two in each step.

Suppose first that \(s=1\). For every nonidentity complement
element, every lift centralizes \(H\) by the abelian coset
formula \eqref{eq:abelian-coset}. The conjugates of the
complement generate \(Q\). Indeed, for a nonidentity
complement element \(a\), the map \(1-a\) is invertible on
the translations, and its values arise as commutators with
\(a\). Because \(H\) is itself abelian, it follows that
\(H\le\Z(G)\). Thus \(K\) is nilpotent, as it is a central
extension of an abelian group, and the preimage of each cyclic
complement is abelian.

A maximum clique in \(K\) can be chosen in \(K\setminus H\).
If \(K\) is nonabelian, its maximum cliques avoid the centre.
If \(K\) is abelian, choose one lift of a nonzero translation.
It follows that \(m+\omega(K)\le N\). The cover of \(K\)
together with the \(m\) abelian complement preimages costs
\(m+T(k)\), where \(k=\omega(K)\).
If \(N\le8\), then \(k\le5\) and the cost is at most
\(N\le T(N)\). For \(N\ge9\) and \(k\le8\), use
\(T(k)\le k+3\) to bound the cost by \(N+3\le \kappa^2\).
If \(k\ge9\), the \(m\ge3\) complements leave \(k\le N-3\),
so \(B(k)\le B(N)/2\). The cost is at most
\(N-8+B(N)/2\le B(N)\), since \(N-8\le2\kappa-6\le2^{\kappa-1}\).
Thus the bound is strict for \(N\ge9\).

Now suppose \(s\ge2\). The complement fibres contribute
at least \(2m\) clique vertices and the translation
fibres contribute at least one, so \(\kappa\ge m\).
In the fibre cover \eqref{eq:affine-kernel-cost}, with
\(a(H)=I=1\), those contributions give a cost
\(1+c+2m\ell\). Each further vertex costs at most
\(D=\max(c,\ell)\) subgroups. Since \(u+v\le N\le2\kappa+2\),
there are at most \(1+2(\kappa-m)\) further vertices. Hence
\[
 a(G)\le 1+c+2m\ell+D+2(m-1)(\kappa-m).
\]
The bounds on \(c,\ell\) give \(1+c+2m\ell+D\le m^2\):
for odd \(m\) use \(c,D,\ell\le(m-1)/2\), and for even
\(m\) use \(c,D\le m-1\) and \(2\ell\le m-2\).
Consequently
\[
 a(G)\le m^2+2(m-1)(\kappa-m)\le \kappa^2.
\]
This proves \(a(G)\le B(N)<T(N)\) when \(\kappa\ge4\).
The only smaller value is \(\kappa=3\), where \(N=7,8\)
and the same bound gives \(a(G)\le9<T(N)\).

\subsection{Nonabelian kernel and cyclic complement}

A narrow fibre supplies a large abelian subgroup of the kernel
and controls the cost of its cosets.

\begin{lemma}[Covering the kernel from one fibre]\label{lem:affine-index}
Let \(H\lhd G\) be nonabelian and let \(gH\) have width \(r\).
Assume \(a(C)\le T(r)\) for each \(C\le H\) with
\(\omega(C)\le r\). There is an abelian subgroup \(A\le H\)
of index \(I\), centralized by an element of \(gH\), with
\[
 a(H)\le rI\le \frac3{16}B(4r),\qquad I\le2^r.
\]
\end{lemma}

\begin{proof}
As in Section~\ref{sec:clique-calculus}, if \(x\in gH\)
centralizes \(A\) then the \(I\) cosets of \(A\) in \(H\)
cost at most \(rI\) in all. We find \(A\) inside a
centralizer in the cover of a maximum clique in \(gH\).
The centralizers in \(H\) have intersection \(D\le\Z(H)\)
and satisfy \([H:D]\le r!\) by \eqref{eq:tomkinson-coset}.

\emph{The uniform construction.}
For \(r\ge4\), if all centralizers have index at least \(r\),
their cosets partition the fibre into commuting sets by
Lemma~\ref{lem:centralizer-partition}. Each centralizer
is then abelian of index \(r\). Otherwise choose
\(C=C_H(x)\) of index \(i\le r-1\).
It satisfies \(\omega(C)\le r\) by Lemma~\ref{lem:coset-calculus}.
If \(C\) is nonabelian, a cover by at most \(T(r)\)
abelian subgroups contains one of index at most \(T(r)-1\):
if the largest has order \(h_0\), their union has size at
most \(1+T(r)(h_0-1)<T(r)h_0\). If \(C\) is abelian, take
\(C\) itself. For \(r\ge5\), this gives \(I\le(r-1)(T(r)-1)\).

\emph{Small widths.}
Only widths three and four need sharper information about
the central quotient.
For \(r=3\), the noncyclic quotient \(H/\Z(H)\) has order
four or six. In the latter case it is \(S_3\), whose
automorphisms are inner. As \(H\) is normal, one can choose
a representative of \(gH\) that centralizes \(H/\Z(H)\)
in \(G/\Z(H)\). Translating the three transpositions and
a \(3\)-cycle by this representative gives a four-point
clique in the image of \(gH\), a contradiction. Thus
\([H:\Z(H)]=4\) and \(D=\Z(H)\). Some centralizer has
index at most three by Lemma~\ref{lem:coset-calculus}.
Its index divides four and is at least two, hence equals two.
Its quotient by \(\Z(H)\) is cyclic, so this centralizer is abelian.

For \(r=4\), the construction above gives \(I\le6\) if \(i\le2\),
and \(I=4\) in the partition case. If \(i=3\), we need an abelian subgroup
of index at most two in \(C\). The same centralizer
intersection gives
\(z(C)\le[C:D]\le4!/3=8\).
If \(C/\Z(C)\) has an element of order at least three,
its cyclic subgroup has index at most two, and its
preimage is abelian. Otherwise the quotient is elementary
abelian of order at most eight. A line suffices when its
order is two or four, and \(C\) itself suffices if it is
abelian. At order eight, two distinct nonzero classes
must have commuting lifts, since \(\omega(C)\le4\).
They span a plane with abelian preimage of index two.
In all cases \(I\le6\).

\emph{The resulting bounds.}
For \(r=3,4\), the index bounds give \(rI\le(3/16)B(4r)\).
For \(r\ge5\), the profile gives \(T(r)-1\le2^{r-3}\), so
\[
 \frac{rI}{B(4r)}\le r(r-1)2^{-r-2}\le \frac5{32}.
\]
The last expression decreases from \(r=5\) since
\((r+1)/(2(r-1))<1\). The same index bounds give \(I\le2^r\):
for \(r\ge5\), both \(r-1\) and \(T(r)-1\) are at most
\(2^{r/2}\).
\end{proof}

Take \(r=\min(s,u)\).
Both widths are at least three,
and \(N\ge v+u\ge(m+1)r\) ensures \(r<N\).
The induction hypothesis gives the cover bound required
in Lemma~\ref{lem:affine-index}. With its subgroup \(A\),
the fibre estimate \eqref{eq:affine-kernel-cost} becomes
\[
 a(G)\le I(r+cu+\ell v).
\]
The translation fibres contribute at least \(r\) clique
vertices, and the complements at least \(mr\). At these widths,
the cover bound is \(Ir(1+c+m\ell)\), including the
\(Ir\) subgroups for the kernel. The bounds on
\(c,\ell\) give \(1+c+m\ell\le(m^2+1)/2\) for either
parity of \(m\). Put \(M=(m^2+1)/2\).

For fixed \(m\), each vertex beyond \((m+1)r\) increases
the cover bound by at most \(I(m-1)\), since \(c,\ell\le m-1\).
As \(rM\ge3(m-1)\), the linear bound
\[
 IrM+I(m-1)\bigl(N-(m+1)r\bigr)
\]
grows by a factor at most \(4/3\) per extra vertex.
The lower bound \(2^{N/2-1}\) for \(B(N)\) grows faster,
by a factor \(\sqrt2\). Comparing at the starting size
\((m+1)r\) therefore gives
\[
 \frac{a(G)}{B(N)}
 \le \frac{rI}{B(4r)}\,M2^{-(m-3)r/2}.
\]
Finally, increasing \(m\) by one adds \(r\) vertices to
the starting clique size. The factor \(M\) increases by a
factor of at most two, whereas the exponential denominator
increases by \(2^{r/2}>2\). Thus the last bound is largest at \(m=3\),
where \(M=5\), and is at most \(15/16\) by
Lemma~\ref{lem:affine-index}.

\subsection{Nonabelian complement}\label{sec:nonabelian-complement}

Choose \(1\ne x\in\Z(L)\). Its fixed-point-free action on
\(V\) gives \(C_Q(x)=L\). We enlarge a clique in this
centralizer through the translation fibres.

Take a maximum clique \(D\) in \(L\). For each \(c\in D\),
use the two complement elements \(c,xc\). Distinct pairs
are completely noncommuting modulo \(V\). Coprimeness gives
\[
 V=U_c\oplus W_c,\qquad
 U_c=\ker(1-c),\quad W_c=(1-c)V.
\]
Fixing the \(U_c\) component and varying the \(W_c\) component
gives a clique of size \(|W_c|\) over \(c\). When \(U_c\ne0\),
choose its component nonzero. Commutation with an element
\((w,e)\), where \(e\) commutes with \(c\), requires
\[
 (1-e)v=(1-c)w.
\]
Projection onto \(U_c\) shows that this clique is completely
noncommuting with both fibres over \(xc\) and over \(x\):
both act without nonzero fixed vectors on \(U_c\).

The fixed spaces of \(c\) and \(xc\) intersect trivially.
If both are nonzero, the two cliques have sizes with
product at least \(m\), so together give at least
\(2\sqrt m\) vertices. If one is zero, its full fibre
gives \(m\ge2\sqrt m\) vertices. Use such a pair for
\(m\) vertices if it exists. Otherwise adjoin the full
fibre over \(x\). Each of at least \(|D|-1\) further
pairs contributes \(2\sqrt m\). Thus
\[
 \omega(Q)\ge m+2\bigl(\omega(L)-1\bigr)\sqrt m
 \ge q_0:=m+4\sqrt m.
\]
A nonabelian \(L\) requires \(d\ge2\). It has a nonabelian
Sylow subgroup, so its order is divisible by the cube of
a prime. Since \(p\nmid|L|\), its order divides
\(\prod_{i=1}^{d}(p^i-1)\), the part of \(|\mathrm{GL}_d(p)|\)
coprime to \(p\). For \(m=4,8,16\) this product is \(3,21\)
and \(3^2\cdot5\cdot7\), respectively. None has a prime
cube divisor. The remaining possibilities are \(m=9\) and
\(m\ge25\). The latter gives \(N\ge q_0\ge45\), already
treated in Section~\ref{sec:centralizer-descent}. Thus \(m=9\),
\(N\ge21\), and \(L\) is a \(2\)-group of order at most \(16\),
since \(|\mathrm{GL}_2(3)|=48\) and \(3\nmid|L|\).

The kernel \(H\) is abelian. Otherwise every fibre has width
at least three by \eqref{eq:tomkinson-coset}, and a clique
of size at least \(21\) in \(Q\) would give \(N\ge63\).
Write \(w(q)\) for the width of the fibre over \(q\).
For \(q=(v,c)\in Q\setminus V\), its conjugates by translations
are the \(|W_c|\) elements \((v+t,c)\) with \(t\in W_c\).
Distinct such elements do not commute, since \(U_c\cap W_c=0\).
Their fibres have the same width, so \(|W_c|w(q)\le N\).
Summing over the \(m=|U_c||W_c|\) elements with linear part
\(c\), and then over \(c\ne1\), gives
\[
 \sum_{q\in Q\setminus V} w(q)\le N\sum_{c\in L\setminus\{1\}}|U_c|.
\]
If \(c\ne1\) and \(U_c\ne0\), its eigenvalues are \(1,-1\),
since the action is coprime and the field is \(\F_3\).
Thus \(|U_c|=3\) and \(\det(c)=-1\). There are at most
\(|L|/2\) such elements by the determinant homomorphism.
For the other \(c\ne1\) we have \(|U_c|=1\), so the sum is
at most \(2|L|-1\le31\). Add one fibre for each of the
four lines in \(V\) and include \(H\) itself. By
Lemma~\ref{lem:fibre-calculus},
\[
 a(G)\le1+4N+31N<B(N)\qquad (21\le N\le44),
\]
since \(1+35\cdot22<2^{10}=B(21)\) and
\(1+35\cdot44<2^{11}=B(23)\).
This completes the proof of Proposition~\ref{prop:nonnil-cover}.

\section{Attaining the profile}\label{sec:opt}

\begin{proof}[Proof of Theorem~\ref{thm:main}]
Propositions~\ref{prop:nilpotent-cover} and
\ref{prop:nonnil-cover} give \(a(G)\le T(\omega(G))\).
Since \(T\) is nondecreasing, \(h(n)\le T(n)\).
The groups in Section~\ref{sec:models} attain these
covering bounds, so \(h(n)=T(n)\).

Fix \(n\ge9\) and suppose \(a(G)=T(n)\) with
\(N=\omega(G)\le n\). Then \(T(N)=T(n)\), so \(N\ge9\)
and \(\lfloor(N-1)/2\rfloor=\lfloor(n-1)/2\rfloor\).
By Propositions~\ref{prop:nilpotent-cover} and~\ref{prop:nonnil-cover},
\(G\) has the binary symplectic model of dimension \(N-1\).
Thus \(N\) is odd and this dimension is \(2\lfloor(n-1)/2\rfloor\).
Conversely, the extraspecial model in this dimension
attains the bound by the symplectic spread.

The central-index theorem \cite[Theorems~1.1 and~1.2(v)]{index}
gives \(\zeta(n)=(n-1)^2\) for \(3\le n\le6\),
\(\zeta(7)=\zeta(8)=81\), and
\(\zeta(n)=4^{\lfloor(n-1)/2\rfloor}\) for \(n\ge9\).
Thus \(T(n)=1+\sqrt{\zeta(n)}\) for \(n\ge3\).
For \(n\ge9\), its equality case is the same binary
symplectic model just obtained for the covering problem.
\end{proof}

\paragraph{Extremal profiles and equality for individual groups.}
The identity \(h(n)=1+\sqrt{\zeta(n)}\) compares two maxima
under the same clique bound. The equality cases show what
that bound selects. For \(q=2^r\), the determinant on
\(\F_q^2\) and its binary trace both give central index
\(q^2\) and covering number \(q+1\). The same field lines
give an optimal cover in each model. Their clique numbers,
however, are \(q+1\) and \(2r+1\), respectively. Thus the
clique constraint can distinguish commutator maps that give
the same central index and the same covering number.

For an individual group satisfying \(a(G)=1+\sqrt{z(G)}\),
let \(s(G)\) be the largest order of the image of an abelian
subgroup in \(G/\Z(G)\). In the tight case, \(s(G)=\sqrt{z(G)}\),
and every optimal cover induces a partition of the nonidentity
elements of the central quotient. In the loose case,
\(s(G)>\sqrt{z(G)}\), and every optimal cover has two members
whose images share a nonidentity element. The loose case
occurs even for groups of class two with elementary abelian
central quotient. The open question is which incidence
properties of the abelian images determine whether exactly
\(1+\sqrt{z(G)}\) members are needed in this regime.

\begin{samepage}
\section*{Statements and declarations}

\paragraph{Funding.}
The author received no specific funding for this work.

\paragraph{Competing interests.}
The author declares no competing interests.

\paragraph{Data availability.}
No datasets were generated or analysed during this study.

\paragraph{Author contributions.}
The sole author carried out the research and wrote the manuscript.

\end{samepage}

\section*{Declaration of generative AI and AI-assisted technologies}
OpenAI ChatGPT assisted with mathematical discussion,
reference checks and \LaTeX{} editing. The author takes responsibility
for the mathematical arguments, references and final text.

\end{document}